\documentclass[12pt]{amsart}
\usepackage{amssymb, amsmath, amscd}

\usepackage{mathptmx}

\newcommand{\mm}{\mathfrak m}

\newcommand{\Z}{\mathbb{Z}}
\newcommand{\R}{\mathbb{R}}
\newcommand{\N}{\mathbb{N}}

\newcommand{\Ac}{\mathcal{A}}
\newcommand{\Bc}{\mathcal{B}}
\newcommand{\Cc}{\mathcal{C}}
\newcommand{\Dc}{\mathcal{D}}

\newcommand{\Fc}{\mathcal{F}}

\newcommand{\Pc}{\mathcal{P}}
\newcommand{\Sc}{\mathcal{S}}

\newcommand{\Tc}{\mathcal{T}}

\DeclareMathOperator{\Ext}{Ext}

\DeclareMathOperator{\gr}{gr}

\DeclareMathOperator{\Hom}{Hom}

\DeclareMathOperator{\id}{id}
\DeclareMathOperator{\Image}{Im}

\DeclareMathOperator{\Ker}{Ker}

\DeclareMathOperator{\rank}{rank}

\DeclareMathOperator{\Tor}{Tor}

\DeclareMathOperator{\tensor}{\otimes}

\DeclareMathOperator{\pnt}{\raise 0.5mm \hbox{\large\bf.}}

\DeclareMathOperator{\lmod}{\!-Mod}
\DeclareMathOperator{\rmod}{Mod-\!}

\DeclareMathOperator{\chcx}{Ch}
\DeclareMathOperator{\lk}{lk}
\DeclareMathOperator{\op}{op}

\DeclareMathOperator{\colim}{\mathrm{colim}}

\newcommand{\xto}{\xrightarrow}

\newtheorem{thm}{\bf Theorem}[section]
\newtheorem{lem}[thm]{\bf Lemma}
\newtheorem{cor}[thm]{\bf Corollary}
\newtheorem{prop}[thm]{\bf Proposition}

\theoremstyle{definition}
\newtheorem{defn}[thm]{\bf Definition}

\newtheorem{rem}[thm]{\bf Remark}
\newtheorem{ex}[thm]{\bf Example}

\title{Cohomology of partially ordered sets and local cohomology of section rings}

\author{Morten Brun}
\address{Department of Mathematics, University of Bergen, 5008 Bergen, Norway}
\email{morten.brun@mi.uib.no}

\author{Winfried Bruns}
\address{FB Mathematik/Informatik, Universit\"at Osnabr\"uck, 49069 Osnabr\"uck, Germany}
\email{winfried@mathematik.uni-osnabrueck.de}

\author{Tim R\"omer}
\address{FB Mathematik/Informatik, Universit\"at Osnabr\"uck, 49069 Osnabr\"uck, Germany}
\email{troemer@mathematik.uni-osnabrueck.de}

\begin{document}

\begin{abstract}
We study local cohomology of rings of global sections of sheafs on
the Alexandrov space of a partially ordered set. We give a criterion
for a splitting of the local cohomology groups into summands
determined by the cohomology of the poset and the local cohomology
of the stalks. The face ring of a rational pointed fan can be
considered as the ring of global sections of a flasque sheaf on the
face poset of the fan. Thus we obtain a decomposition of the local
cohomology of such face rings. Since the Stanley-Reisner ring of a
simplicial complex is the face ring of a rational pointed fan, our
main result can be interpreted as a generalization of Hochster's
decomposition of local cohomology of Stanley-Reisner rings.
\end{abstract}

\maketitle

\section{Introduction}
In this paper we study local cohomology
of rings of global sections of sheafs on the Alexandrov space of a partially ordered set.
Before we introduce the concepts needed to state our main results,
we describe some consequences.

Recall that a {\em simplicial complex} $\Delta$ on a finite vertex set $V$ is
a collection $\Delta$ of subsets of $V$ closed under inclusion, that
is, if $F \subseteq G$ and $G \in \Delta$, then $F \in \Delta$.
A rational pointed fan $\Sigma$ in $\R^d$ is
a finite collection of rational pointed cones in $\R^d$
such that
for $C' \subseteq C$ with $C\in \Sigma$
we have that $C'$ is a face of $C$ if and only if $C' \in \Sigma$,
and
such that if $C,C' \in \Sigma$, then $C\cap C'$ is a common face of $C$ and $C'$.
The face poset $P(\Sigma)$ of $\Sigma$ is the partially ordered
set of faces of $\Sigma$ ordered by inclusion.
The {\em face ring $K[\Sigma]$} of $\Sigma$ over a field $K$ is
defined as follows:
As a $K$-vectors space $K[\Sigma]$ has one basis element $x^a$ for
each $a$ in the intersection of $\Z^d$ and the union of the faces of
$\Sigma$. Multiplication in $K[\Sigma]$ is defined by
\begin{displaymath}
  x^a x^b =
  \begin{cases}
    x^{a+b} & \text{if $a$ and $b$ are elements of a common face of
    $\Sigma$,} \\
    0 & \text{otherwise.}
  \end{cases}
\end{displaymath}

A simplicial complex $\Delta$ on the vertex set
$V = \{1,\dots,d-1\}$ has an associated rational fan $\Sigma(\Delta)$
with one cone $C(F)$ for each face $F$ of $\Delta$. The cone
$C(F)$ is equal to the set of
vectors
in $\R_{\ge 0}^d = \{(x_1,\dots,x_d) \colon x_1,\dots,x_d \ge 0\}$
with $x_d =\sum_{i=1}^{d-1} x_i$ and $x_i = 0$ for $i \in V \setminus F$.
The face ring $K[\Sigma(\Delta)]$ is
called the {\em Stanley-Reisner ring $K[\Delta]$} of $\Delta$.

Stanley observed in \cite[Lemma 4.6]{ST87} that the face ring
$K[\Sigma]$ of a rational pointed fan $\Sigma$ is a Cohen-Macaulay
ring if the Stanley-Reisner ring $K[\Delta(P(\Sigma))]$ of the order
complex of the face poset of $\Sigma$ is a Cohen-Macaulay ring.
Stanley reduces the proof of this observation to a theorem of
Yuzvinsky \cite[Theorem 6.4]{YU87}. Theorem \ref{corgenhochster-1}
below generalizes Yuzvinsky's theorem, and Stanley's
  observation is a direct consequence of it.

We consider every poset $P$ as a topological space with the {\em Alexandrov
  topology}, that is, the
topology where the open sets are the  lower subsets (also called order ideals) of $P$.
If $\Tc$ is a sheaf of $K$-algebras on $P$, then for every $x \in P$ there is
a restriction homomorphism $H^0(P,\Tc) \to \Tc_x$ from the zeroth
cohomology group
of $P$ with coefficients in $\Tc$ to the
stalk $\Tc_x$ of $\Tc$ at $x$. Given an ideal $I$ in
a commutative ring $R$ and an $R$-module $M$ we denote the local
cohomology groups of $M$ by $H^i_I(M)$ for $i \ge 0$.
The following is our  decomposition of
local cohomology.
\begin{thm}
\label{genhochster-1}
Let $K$ be a field, let $\Tc$ be a sheaf of $K$-algebras on a finite poset $P$ and
let $I$ be an ideal of the
zeroth cohomology group
$H^0(P,\Tc)$ of $P$ with coefficients in $\Tc$.
For $x \in P$ we let
$d_x$ denote the Krull dimension of the stalk $\Tc_x$ of $\Tc$ at $x$
and we
assume that:
\begin{enumerate}
\item
$H^0(P,\Tc)$ is a Noetherian ring
and
$H^i(P,\Tc) = 0$ \ for every $i > 0$,
\item
$H^i_I(\Tc_x) = 0$ for every $x \in P$ and every $i \ne d_x$,
\item
if $x < y$ \ in $P$ then $d_x<d_y$.
\end{enumerate}
Then
there is an isomorphism
$$
H^i_I(H^0(P,\Tc))
\cong
\bigoplus_{x \in P}
\widetilde H^{i - d_x -1}((x,1_{\widehat P});K)
\tensor_K
H^{d_x}_I(\Tc_x)
$$
of $K$-modules, where $\widetilde H^{i - d_x -1}((x,1_{\widehat
  P});K)$ denotes the reduced cohomology of the partially ordered set
$(x,1_{\widehat P}) = \{y \in P \colon x < y\}$ with coefficients in
$K$. If $\Tc$ is a sheaf of $\Z^d$-graded $K$-algebras, the above
isomorphism is an isomorphism of $\Z^d$-graded $K$-modules.
\end{thm}
The zeroth cohomology ring $H^0(P,\Tc)$ is naturally identified with the ring of
global sections of $\Tc$. Under this name it was studied by
Yuzvinsky \cite{YU87, YU89} and
Caijun \cite{CA97}. The following
immediate corollary of Theorem \ref{genhochster-1}
generalizes the results
\cite[Theorem 2.4]{CA97} and \cite[Theorem  6.4]{YU87} of Cajun
and Yuzvinsky.
\begin{thm}
\label{corgenhochster-1}
  Suppose in the situation of Theorem \ref{genhochster-1} that there
  exists a unique graded maximal ideal
  $\mm$ in $H^0(P,\Tc)$ which is maximal considered as an ideal of $H^0(P,\Tc)$.
  If
  the assumptions of Theorem \ref{genhochster-1} are satisfied by the
  ideal $\mm$, then
  the ring
  $H^0(P, \Tc)$ is Cohen-Macaulay if and only if there exists
  a number $n$ such that
  the reduced cohomology $\widetilde H^*(
  \Delta((x,1_{\widehat P})),K)$ of the simplicial
  complex $\Delta((x,1_{\widehat P}))$ associated to
  the poset $(x,1_{\widehat P})$ is concentrated in  degree
  $(n-d_x-1)$ for every $x \in P$ with $\Tc_x \ne 0$.
\end{thm}

Let us return to the situation where $P = P(\Sigma)$ is the face poset
of a rational pointed fan $\Sigma$ in $\R^d$. There is a unique
flasque sheaf $\Tc$ of $\Z^d$-graded $K$-algebras associated to
$\Sigma$ with stalks $\Tc_C = K[C]$ given by the monoid algebras on
the cones of $\Sigma$, with a natural isomorphism $K[\Sigma] \cong
H^0(P,\Tc)$ between the face ring of $\Sigma$ and the ring of global
sections of $\Tc$ and with restriction homomorphisms $H^0(P,\Tc) \cong
K[\Sigma] \to K[C] \cong \Tc_C$ acting by the identity on $x^a$ if $a
\in C$ and  taking $x^a$ to zero otherwise
(see \cite[Theorem 4.7]{BRRO04}).
Note that
the face ring $K[\Sigma]$
is Noetherian, $\Z^d$-graded and has a unique maximal $\Z^d$-graded ideal $\mm$.
The $K$-algebras $K[C]$ are normal
and thus Cohen-Macaulay of Krull dimension $d_C = \dim(C)$.
Hence the following result is a direct consequence of Theorem \ref{genhochster-1}
since
every flasque sheaf of rings on a poset satisfies assumption (i).

\begin{thm}
Let $\Sigma$ be a rational pointed fan in $\R^d$ with face poset $P$,
$K$ be a field and
$\mm$ be the graded maximal ideal of the face ring $K[\Sigma]$.
Then there is an isomorphism
$H^i_\mm(K[\Sigma]) \cong \bigoplus_{C \in P}
\widetilde{H}^{i-\dim(C)-1}((C,1_{\widehat P});K)
\otimes_K
{H^{\dim(C)}_\mm(K[C])}
$
of $\Z^d$-graded $K$-modules.
\end{thm}
If $\Sigma = \Sigma(\Delta)$ for
a simplicial complex $\Delta$, then the
posets $(C,1_{\widehat P}) = (C(F),1_{\widehat P})$ in the above formula
are isomorphic to the face posets $P(\lk_{\Delta} F) \setminus \emptyset$ of the
links in $\Delta$. Thus the above theorem generalizes Hochster's
decomposition
of local cohomology of Stanley-Reisner rings
(see \cite{BRHE98}).
Just like Reisner's topological characterization of
the Cohen-Macaulay property  of Stanley-Reisner rings is a consequence of Hochster's decomposition
of local cohomology of Stanley-Reisner rings, the observation
\cite[Lemma 4.6]{ST87}
of Stanley mentioned above is a consequence of Theorem \ref{corgenhochster-1}.

Sheaves on a poset
$P$ can be described in a different way.
More precisely,
let $R$ be a commutative ring,
then a sheaf $\Tc$ of $R$-algebras
on $P$ is described by a unique collection $(T_x)_{x \in P}$ of
$R$-algebras
and homomorphisms
$T_{xy} \colon T_y \to T_x$ for $x \le y$ in $P$ with
the property that $T_{xx}$ is the identity on $T_x$
and that $T_{xy}\circ T_{yz} = T_{xz}$ for every $x \le y \le z$ in $P$.
Moreover, every such collection describes a unique sheaf on $P$,
and we have that $H^0(P,\Tc)$ is the (inverse) limit $\lim T_x$.
In the body of this paper we call $T = (T_x,T_{xy})$
an {\em $RP$-algebra}, and we work with $RP$-algebras instead of with sheaves.
One reason for this change of perspective is
that homological algebra of $RP$-algebras
is more accessible than sheaf cohomology.
In fact,
Theorem \ref{genhochster-1} is a consequence of general homological arguments.

Algebras of  type $\lim T$ for an $RP$-algebra $T$
appear at many places
in commutative algebra and combinatorics.
For example Bruns and Gubeladze  studied such algebras in a series of papers
\cite{BRGU99, BRGU01}.
Brun and R\"omer considered
the relationship between
initial ideals
of the defining ideal
of the face ring of a rational fan
and
subdivisions of that fan
in \cite{BRRO04}.

This paper is organized as follows.
In Section \ref{prere} we recall definitions and notations
related to posets and abstract simplicial complexes.
In Section \ref{examplesrp} we introduce $RP$-algebras and give
examples from commutative algebra and combinatorics.
In Section \ref{mainresults} we study the local cohomology
of limits of $RP$-algebras.
In particular, we prove Theorem \ref{genhochster-1}.
In Section \ref{application} we present further applications.

The proof of Theorem \ref{genhochster-1} uses
homological algebra over a poset $P$.
We recollect in Sections \ref{kpmodules} and \ref{homologicalalg}
some notation and basic facts needed in the proof. Apart from
Proposition \ref{sheafpmod}
the results in Sections \ref{kpmodules} and \ref{homologicalalg}  are well-known.
We have included these sections as a
bridge between the questions
considered
in this paper and
the literature on homological algebra of functor categories.
They are also intended as a
soft introduction to the language of
functor categories for the reader who is not very
familiar with category theory.
More detailed accounts of the concepts introduced here can be found in
the papers of Baues-Wirsching \cite{BW} and of Mitchell \cite{Mitchell}.
Sections \ref{locorp} is the technical heart of the paper.
Here we use the theory collected in Sections \ref{kpmodules} and \ref{homologicalalg}
to study the local cohomology of $RP$-modules.
\section{Prerequisites}
\label{prere}
In this paper $P = (P,\le)$
always denotes a {\em partially ordered set} (poset for short).
Given $P$,
the opposite poset $P^{\op} = (P,\preceq)$ has the same underlying set as
$P$ and the opposite partial order $\preceq$, that is,
$y \preceq x$ if and only if $x \le y$.
We write $x < y$ if $x \neq y$ and $x \le y$.
We also consider a poset as a category ``with morphisms pointing down'', that
is, for $x,y \in P$ there is a unique arrow $y \to x$ if and only if
$x \le y$.
If $P$ contains a unique maximal element, then this element
is called the {\em initial element} of $P$, and it will be denoted $1_P$.
Analogously a unique minimal element of $P$ is called a
{\em terminal element} of $P$ and it is denoted $0_P$.
The poset $\widehat P = (\widehat P,\le)$ associated with $P$ has
underlying set $\widehat
P = P \cup \{0_{\widehat P},1_{\widehat P}\}$ obtained by adding a
terminal element $0_{\widehat P}$ and an initial
element $1_{\widehat P}$ to $P$ (in spite of a terminal or an initial element
which may already exist in $P$).
The closed interval $[x,y]$ of elements between $x$ and $y$ in $P$
is the set $[x,y] = \{z \in P \colon x \le z \le y\}$ considered as
a sub-poset of $P$.
The half-open interval $[x,y)$,
the half-open interval $(x,y]$
and the open interval $(x,y)$ are
described similarly.
Note that
if $x \in \widehat P$ and $y \in P$, then $(x,y]$ is a sub-poset of $P$.
A finite poset $P$ is called a {\em graded} poset,
if all maximal chains (i.e. totally ordered subsets) of $P$ have the same length $\rank(P)$.
In this situation it is possible to define a unique rank function on $P$
such that for $x \in P$ we have that $\rank(x)$ is the common
length of maximal chains in $P$ ending at $x$.

The {\em Alexandrov topology} \cite{Al56} on a poset $P$ is the
topology where the
open subsets are the lower subsets (also called order ideals),
that is, the subsets
$U$ such that $y \in U$ and $x \le y$ implies $x \in U$.
The subsets of the form $(0_{\widehat P},x]$ form a basis for this
topology.
(Sometimes, e.g. in \cite{YU87}, this is called the
{\em order topology} on $P^{\op}$.)

The poset $P$ is {\em locally finite} if every closed interval of the
form $[x,y]$ for $x,y \in P$ is finite and it is
{\em topologically finite} if every interval of the form $(0_{\widehat P},x]$ for $x
  \in P$ is
finite.
Note that locally finite does not imply topologically finite but the
converse is true
and that $P$ is topologically finite if and only if
every element of $P$ has a finite neighborhood in
the Alexandrov topology.

If a simplicial complex $\Delta$ on a finite vertex set $V$
is non-empty, then the empty set is
a terminal element in the {\em face poset}, that is, in the partially
ordered set $P(\Delta) =
(\Delta, \subseteq)$ of elements in $\Delta$ ordered by inclusion.
The elements $F$ of $\Delta$ are called {\em faces}. If $F$ contains $d+1$
vertices, that is, $d+1$ elements of $V$, then $F$ is called a
{\em $d$-dimensional face}, and
we write $\dim F = d$. The empty set is a face of
dimension $-1$.
The {\em dimension} $\dim \Delta $
is the supremum of the dimensions of
the faces of $\Delta$.

Above we constructed a poset $P(\Delta)$ associated to every
simplicial complex $\Delta$.
Conversely, the {\em order complex}
$\Delta(P)$ of $P$ is the simplicial complex on the vertex set $P$
consisting of the chains in $P$ ordered by inclusion.

Recall that if $F$ is a face of a simplicial complex $\Delta$ then the
{\em link} $\lk_\Delta F$ of $F$ in $\Delta$ is
the simplicial complex
$\lk_\Delta F=\{ G \setminus F : G \in
\Delta \text{ and } F \subseteq G \}$.
It is easy to see that
the correspondence
$G \setminus F \mapsto G = (G \setminus F) \cup F$
defines an order-preserving bijection
$P(\lk_{\Delta} F) \xto \cong [F,1_{\widehat {P(\Delta)}})$.
Thus the barycentric subdivision
$\Delta(P(\lk_{\Delta} F) \setminus \emptyset )$
of $\lk_{\Delta} F$
has the same reduced simplicial (co-) homology
as $\Delta((F,1_{\widehat {P(\Delta)}}))$.
This fact will be used
several times in this paper.

For more details on simplicial complexes
and posets
see for example the corresponding chapters in the books of
Bruns-Herzog \cite{BRHE98} and Stanley \cite{ST96,ST99}.

\section{Examples of $RP$-Algebras}
\label{examplesrp}
Fix a commutative ring $R$ and a poset $(P,\le)$.
An {\em $RP$-algebra} $T$ is a system $(T_x)_{x \in P}$ of $R$-algebras and
homomorphisms $T_{xy} \colon T_y \to T_x$ for $x \le y$ in $P$ with
the property that $T_{xx}$ is the identity on $T_x$ and that $T_{xy}
\circ T_{yz} = T_{xz}$ for every $x \le y \le z$ in $P$.
The (inverse) {\em limit} of $T$ is the subring $\lim T$ of $\prod_{x \in P} T_x$
consisting of sequences $r=(r_x)_{x \in P}$ with the property that
$T_{xy}(r_y) = r_x$ for every $x \le y$ in $P$. In particular, an
$RP$-algebra is an $RP$-module in the sense explained in Section
\ref{kpmodules} and we can apply the theory developed in that section.
The $RP$-algebra $T$ is called {\em cyclic} if the homomorphism
$R \to T_x$ is surjective for every $x \in P$.
We call
$T$ an $RP$-algebra of $\Z^d$-graded $RP$-algebras if
the $R$-algebras $T_x$ are $\Z^d$-graded,
and the homomorphisms $T_{xy}\colon T_y \to T_x$ and $R \to T_x$ are homogeneous of degree zero for $x \le y$.
In this case $\lim T$ is a $\Z^d$-graded $R$-algebra.
Similarly, we call $T$ an $RP$-algebra of Cohen-Macaulay rings,
if all $T_x$ are Cohen-Macaulay rings.

The following examples motivate
the study of $RP$-algebras in commutative algebra and
algebraic combinatorics.
For a field $K$ and a set $F$ we let $K[F]=K[x_i : i \in F]$
be the polynomial ring
with one indeterminate for each $i \in F$.

\begin{ex}[Stanley-Reisner ring]
\label{stanleyreisnerring}
Let $K$ be a field,
$\Delta$ be a simplicial complex on the vertex set $V=\{1,\dots,d\}$ and
$P=P(\Delta)$.
For $F \in P$ define $T_F=K[F]$.
For $G \subseteq F$, we define $T_{GF}\colon K[F] \to K[G]$
to be the natural projection.
If we let $R = K[\cup_{F \in P} F]$,
then $T$ is a cyclic $RP$-algebra
and
$$
\lim T \cong R/ I_{\Delta},
$$
where
$I_{\Delta}$ is generated by all squarefree monomials
$\prod_{i \in G}x_i$ for $G \subseteq V, G \not\in \Delta$.
Hence $\lim T$ is the usual {\em Stanley-Reisner ring} in this
case.

Note that the polynomial algebra $R =
K[x_1,\dots,x_d]$ is $\Z^d$-graded and $T$
is a cyclic $RP$-algebra of $\Z^d$-graded $R$-algebras.
\end{ex}

\begin{ex}[Toric face rings]
\label{toricrings}
We consider a rational pointed fan $\Sigma$ in $\R^d$,
that is,
a collection of rational pointed cones in $\R^d$
such that
for $C' \subseteq C$ with $C\in \Sigma$
we have that $C'$ is a face of $C$ if and only if $C' \in \Sigma$,
and if $C,C' \in \Sigma$, then $C\cap C'$ is a common face of $C$ and $C'$.
Let $P$ be the face poset of $\Sigma$ ordered by inclusion.
For a rational pointed cone $C \in P$
we let $T_C$ be the
monoid ring $K[C\cap \Z^d]$ over a field $K$.
The homomorphisms $T_{C'C}\colon T_{C} \to T_{C'}$
are induced by the natural face projection.
Then $T$ is a $\Z^d$-graded $KP$-algebra.
For a suitable polynomial ring $R$ over $K$
the homomorphisms $R \to T_C$ are all surjective and $\Z^d$-graded.
Then $T$ is a cyclic $RP$-algebra of $\Z^d$-graded $R$-algebras.
In any case, $\lim T$ is the {\em toric face ring} of $\Sigma$.
This examples goes back to a construction of Stanley in \cite{ST87}.
It was generalized, and these rings were studied
by Brun-R\"omer in \cite{BRRO04}.
The limits $\lim T$
of such algebras were intensively studied by Bruns-Gubeladze
(see \cite{BRGU99, BRGU01}).

Note that Stanley-Reisner rings of simplicial complexes are
covered by the this example: Let $\Delta$ be a simplicial
complex on the vertex set $V = \{1,\dots,d-1\}$. To a subset
$F$ of $V$ we associate the pointed cone $C_F$ in $\R^d$ generated by
the set of
elements of the form $e_i + e_d$ for $i \in F$.
Here $e_i$ denotes the standard basis vector $e_i = (0,\dots,0,1,0,\dots,0)$ of
$\R^d$. If $\Sigma$ denotes the
fan in $\R^d$ consisting of the cones $C_F$ for $F \in \Delta$, then
the face posets of $\Delta$ and of $\Sigma$ are isomorphic and the
$RP$-algebras of the Examples \ref{stanleyreisnerring} and \ref{toricrings}
correspond to each other via this isomorphism.
\end{ex}

Borrowing notation from the theory of sheaves
we call an $RP$-algebra $T$ {\em flasque}
if $\lim T|_U \to \lim T|_V$
is surjective for all open sets
$V \subseteq U$ of the poset $P$.
In the third example we present
a general construction to produce flasque $RP$-algebras.
(See also \cite{YU87}.)
\begin{ex}
\label{ex3}
Let $R$ be a commutative ring and $\Dc$ a distributive lattice of
ideals in $R$ (with respect to sum and intersection). Moreover, let $P$
be a finite subset of $\Dc$ such that $I+J\in P$ for all $I,J\in P$,
and consider $P$ as a poset with $I \geq J$ if $I \subseteq J$.
Let $T$ be the $RP$-algebra given by $T_I=R/I$ for all $I\in
P$ and $T_{JI}$ the natural epimorphism $R/I\to R/J$.
For $I_1,\dots,I_n$ in $P$ let $U$ be the smallest open subset of $P$
containing $I_1,\dots,I_n$. We claim that
$$
\lim T|_U=R/(I_1\cap\dots\cap I_n).
$$
In fact, clearly $\lim T|_U$ is the kernel of the map
$$
\Phi: R/I_1\times\dots\times R/I_n\to \prod_{i<j} R/(I_i+I_j)
$$
where $\Phi(\bar a_1,\dots,\bar a_n)=(\bar a_i- \bar a_j:i<j)$ and
$\bar{\phantom{x}}$ denotes the residue class with respect to the
appropriate ideal.
For $n=2$ the claim is proved by the classical exact sequence
$$
0\to R/(I_1\cap I_2)\to R/I_1\times R/I_2\to R/(I_1+I_2)\to 0.
$$
By induction we can assume that an element in the kernel of $\Phi$
has the form $(\bar a,\dots,\bar a,\bar b)$, and it remains to show
that the sequence
$$
0\to R/\bigl((I_1\cap\dots\cap I_{n-1})\cap I_n\bigr) \to
\bigl(R/I_1\cap\dots\cap I_{n-1}\bigr)\times R/I_n
\xrightarrow{\Phi'} \prod_{i=1}^{n-1} R/(I_i+I_n)
$$
is exact. By the case $n=2$, it is enough that the target of $\Phi'$
can be replaced by $R/(I_1\cap\dots\cap I_{n-1})+I_n$, and this
follows immediately from distributivity.
Finally, it is now immediate that $T$ is flasque.

Distributive lattices of
ideals in rings (with respect to sum and intersection)
appear naturally in commutative algebra.
For example let $S=K[x_1,\dots,x_n]$ be the polynomial ring.
Assume that we have fixed a standard $K$-basis $\Sc$ of $S$
and $\Dc$ is the set of ideals which have a subset of $\Dc$ as
a $K$-basis. Then $\Dc$ is a distributive lattice.
If $\Sc$ is the set of usual monomials in $S$, then
$\Dc$ is the set of monomial ideals.
Choosing a finite set $\Pc$ of monomial prime ideals
closed under summation of ideals gives back the example of Stanley-Reisner rings.
\end{ex}

Example \ref{ex3} is in fact quite general.
First let us note that if
$T$ is a flasque $RP$-algebra on a finite poset $P$
with the property that the
homomorphism $R \to \lim T$ is surjective, then the kernel of this
homomorphism is given by the intersection of the kernels $I_x$ of
the compositions $R \to \lim T \to T_x$ for $x \in P$. Indeed, the
homomorphism $R/\cap_{x \in P} I_x \to \lim T$ is surjective, and
it is injective since the composition $R/\cap_{x \in P} I_x \to \lim
T \to \prod_{x \in P} T_x$ is injective. Next note that if $Q$
denotes the poset consisting of nonzero sums of ideals of the form
$I_x$ for $x \in P$, with order given by reverse inclusion, there is
an order-preserving map $P \to Q$ taking $x\in P$ to $I_x\in
Q$. Since $T_x \cong R/I_x$, there is a canonical homomorphism
$\lim_{I \in Q} R/I \to \lim T$. Since the kernel of $R \to \lim_{I
\in Q} R/I$ is equal to the kernel of $R \to \lim T$, we have
an isomorphism $\lim_{I \in Q} R/I \cong \lim T$. Finally, let us
note that the map $P \to Q$ is injective if and only if $x <y$ in
$P$ implies that the homomorphism $T_y \to T_x$ has a nonzero
kernel.
For further examples of $RP$-algebras see \cite{BRRO05}.


\section{Decomposition Results}
\label{mainresults}
This section contains our main results
on the local cohomology of the limit of an $RP$-algebra $T$.
The following theorem is our most general
decomposition of
local cohomology groups.
The proof will be given in
Section \ref{locorp}.

\begin{thm}
\label{genhochster02}
Let $R$ be a Noetherian commutative algebra over a field $K$
and $I \subset R$  an ideal.
Let $P$ be
a finite poset and let $T$ be an $RP$-algebra. For $x \in P$ we let
$d_x$ denote the Krull dimension of $T_x$ and we
assume that:
\begin{enumerate}
\item
$H^i_I(T_x) = 0$ for every $x \in P$ and every $i \ne d_x$,
\item
$\Ext^i_{RP}(R,T) = 0$ \ for every $i > 0$,
\item
if $x < y$ \ in $P$ then $d_x<d_y$.
\end{enumerate}
Then
there is an isomorphism
$$
H^i_I(\lim T)
\cong
\bigoplus_{x \in P}
\widetilde H^{i - d_x -1}((x,1_{\widehat P});K)
\tensor_K
H^{d_x}_I(T_x)
$$
of $K$-modules.

If
$R$ is a Noetherian commutative $\Z^d$-graded algebra over a field $K$,
$I \subset R$ is a graded ideal
and $T$ is an $RP$-algebra of $\Z^d$-graded $R$-algebras,
then the above isomorphism is an isomorphism of $\Z^d$-graded
$K$-modules.
\end{thm}

Recall that a graded maximal ideal in a
$\Z^d$-graded ring
$R$
is a graded ideal $\mm$ that is maximal among the proper
graded ideals in $R$, and
that $R$ is graded local if it
contains a unique graded maximal
ideal.
In the next section we give applications
of the following (graded) version of
Theorem \ref{genhochster02}.
\begin{cor}
\label{genhochster}
Let $R$ be a Noetherian $\Z^d$-graded local commutative algebra
over a field $K$ with
graded maximal ideal $\mm$
which is maximal considered as an ideal of $R$.
Let $P$ be
a finite poset and let $T$ be a cyclic $RP$-algebra of
$\Z^d$-graded
$R$-algebras.
For $x \in P$ we let
$d_x$ denote the Krull dimension of $T_x$ and we
assume that:
\begin{enumerate}
\item
$T_x$ is a Cohen-Macaulay ring for every $x \in P$,
\item
$\Ext^i_{RP}(R,T) = 0$ \ for every $i > 0$,
\item
if $x < y$ \ in $P$ then $d_x<d_y$.
\end{enumerate}
Then
there is an isomorphism
$$
H^i_\mm(\lim T)
\cong
\bigoplus_{x \in P}
\widetilde H^{i - d_x -1}((x,1_{\widehat P});K)
\tensor_K
H^{d_x}_\mm(T_x)
$$
of $\Z^d$-graded $K$-modules.
\end{cor}
\begin{proof}
  We first note that under the given assumptions $T_x$ is
  Cohen-Macaulay if and only if $(T_x)_\mm$ is Cohen-Macaulay.
  (See \cite[Exercise 2.1.27]{BRHE98}.) By Grothendieck's
  Vanishing Theorem on local cohomology modules
  (see \cite[Theorem 3.5.7]{BRHE98}) this is the case if and only if
  condition (i) of
  Theorem \ref{genhochster02} is satisfied. The other
  conditions of Theorem \ref{genhochster02} are part of our assumptions.
\end{proof}


\section{Applications}
\label{application}
In this section we explain how the results of Section \ref{mainresults}
generalize the
Hochster formulas
for local cohomology.
Let $P$ be a poset and let $R$ be a commutative ring.
The following lemma is a consequence of Lemma \ref{flasquecrit0}:
\begin{lem}
\label{flasquecrit}
Let $R$ be a commutative ring and let $P$ be a poset.
If $T$ is a flasque $RP$-algebra, then
$\Ext_{RP}^i(R,T)=0 \text{ for all } i>0.$
\end{lem}
Using Corollary \ref{genhochster} it is possible to
study several properties of the ring $\lim T$.
The next result generalizes results of
Yuzvinsky \cite[Theorem 6.4]{YU87}
and Caijun \cite[Theorem 2.4]{CA97}
in the graded situation.
\begin{cor}
\label{cmring}
Let $R$ be a $\Z^d$-graded local Noetherian algebra
over a field $K$
with unique graded maximal ideal $\mm$
which is maximal considered as an ideal of $R$.
Assume that $T$ is a flasque cyclic $RP$-algebra of
$\Z^d$-graded  Cohen-Macaulay $R$-algebras such that $x < y$ in $P$ implies
$d_x<d_y$.
The following statements are equivalent:
\begin{enumerate}
\item
$\lim T$ is a Cohen-Macaulay ring,
\item
$
\widetilde{H}^{p}\bigl((x,1_{\widehat{P}});K\bigr)=0$
for $x \in P$ and
$p \neq \dim (\lim T) - d_x-1.$
\end{enumerate}
\end{cor}
\begin{proof}
Recall that by standard arguments
$\lim T$ is a Cohen-Macaulay ring if and only if
$H_\mm^i(\lim T)=0$ for $i\neq \dim (\lim T)$.
The equivalence of (i) and (ii) is a direct consequence of
Corollary \ref{genhochster} and
Lemma \ref{flasquecrit}.
\end{proof}

We need the following result.
\begin{lem}
\label{helperflasque}
Let $T$ be a $KP$-algebra. Assume
that there exists $x \in P$ with the following properties:
\begin{enumerate}
\item
$T_y = T_x$ for every $y \in [x,1_{\widehat P})$,
\item
$T_{xy}$ is the identity homomorphism on $T_x$
for every $y \in [x,1_{\widehat P})$,
\item
$T_y = 0$ if $y \notin [x,1_{\widehat P})$.
\end{enumerate}
Then $T$ is a flasque $KP$-algebra.
\end{lem}
\begin{proof}
If $U$ is an open subset of $P$ then $T|_U = 0$ if $x \notin U$, and
otherwise $\lim T|_U \cong T_x$. If $V \subseteq U$ is an inclusion of
open subsets of $P$ then the natural projection $\lim T|_U \to \lim
T|_V$ is isomorphic to the identity on $T_x$ if $x \in V$, and
otherwise $\lim T|_V = 0$.
Thus $T$ is a flasque $KP$-algebra.
\end{proof}

The following result generalizes
an observation of
Yuzvinsky  \cite[Proposition 7.6]{YU87}.
\begin{prop}
\label{flsh_st}
Let $\Sigma$ be a rational pointed fan in $\R^d$ with face lattice
$P$. The
$\Z^d$-graded
$KP$-algebra $T$ of Example \ref{toricrings} is flasque and if $D \subset
C$ in $P$, then the Krull dimension of $T_D$ is strictly less
than the Krull dimension of $T_C$.
In particular, for
every simplicial complex $\Delta$, the $\Z^d$-graded
$KP(\Delta)$-algebra of Example \ref{stanleyreisnerring} is flasque.
\end{prop}
\begin{proof}
$T$ has the decomposition
$T = \bigoplus_{a \in \Z^d} T(a)$,
where $T(a)_C = K$ if $a \in C$ and $T(a)_C = 0$ otherwise.

Let $D = \cap_{C \in P \colon a \in C} C$. Then $T(a)_{C} = 0$
unless $D \subseteq C$, and in this case $T(a)_{C} = K$ and the map
$T(a)_{DC}$ is the identity map on $K$.
It follows from \ref{helperflasque}
that $T$ is a direct sum of flasque $KP$-algebras,
and this implies that $T$ is a flasque $KP$-algebra.

The statement about Krull dimensions
holds
since
the Krull dimension of the monoid ring $K[C \cap \Z^d]$ is equal to the dimension
of the cone $C$.
\end{proof}
Let $\Delta$ be a simplicial complex with face lattice $P$. We have
observed that for every $F \in P$ the posets
$P(\lk_{\Delta} F)$
and $[F,1_{\widehat P})$ are isomorphic.
Considering
the rational pointed fan induced by a simplicial complex as in Example \ref{toricrings},
and using the barycentric subdivision
homeomorphism, the following theorem recovers Hochster's
decomposition of the local cohomology of Stanley-Reisner rings.

\begin{thm}
\label{hochsterformula}
Let $\Sigma$ be a  rational pointed fan in $\R^d$ with face lattice
$P$.
If $K[\Sigma]$ denotes the toric face ring of $\Sigma$ over a field
$K$ and if
$\mm$ denotes the graded maximal ideal of $K[\Sigma]$,
then there is an isomorphism
\begin{displaymath}
  H^i_\mm(K[\Sigma]) \cong \bigoplus_{C \in P}
  \widetilde{H}^{i-\dim C -1}((C,1_{\widehat P});K) \otimes_K {H^{\dim C }_\mm(K[C\cap \Z^d]))}
\end{displaymath}
of $\Z^d$-graded $K$-modules.
In particular,
$K[\Sigma]$ is a Cohen-Macaulay ring if and only if
$
\widetilde{H}^{p}\bigl((C,1_{\widehat{P}});K\bigr)=0$
for $C \in P$ and
$p \neq \dim (K[\Sigma]) - d_C-1,$
\end{thm}
\begin{proof}
Let $T$ be the $KP$-algebra of Example \ref{toricrings}.
We have that
$$
H^i_{\mm} (K[\Sigma]) \cong  H^i_{\mm}(\lim T).
$$
For $C\in P$ the ring $T_C=K[C\cap \Z^d]$ is a normal monoid ring of Krull
dimension $d_C = \dim C$. Thus it is Cohen-Macaulay
(see \cite[Theorem 6.3.5]{BRHE98}).
By Lemma
\ref{flsh_st}
the $KP$-algebra $T$ is flasque and $d_D<d_C$ for $D \subset C$ in $P$.
Since $T_C$ is a Cohen-Macaulay ring for every $C \in P$
it follows from
\ref{genhochster}
that there is an isomorphism
\begin{eqnarray*}
{H^i_\mm(K[\Sigma])}
&\cong&
\bigoplus_{C \in P}
\widetilde{H}^{i - \dim C -1}((C,1_{\widehat{P}});K) \otimes {H^{\dim C}_\mm(K[C\cap \Z^d])}
\end{eqnarray*}
of $\Z^d$-graded $K$-modules for every $i \ge 0$.
\end{proof}

\begin{rem}
\label{reisner}
Reisner's Cohen-Macaulay criterion
for a simplicial complex $\Delta$
states that $K[\Delta]$ is a Cohen-Macaulay ring if and only if
$$
\Tilde{H}^i(\lk_\Delta F;K)\cong
\Tilde{H}_i(\lk_\Delta F;K)=0 \
\text{ for all }\
F\in \Delta\ \text{ and all }\ i<\dim\lk_\Delta F.
$$
This criterion is a direct consequence of
Hochster's
decomposition of the local cohomology of
Stanley-Reisner rings.
In particular, since
$
\widetilde{H}^{p}(\lk_{\Delta} F;K)
\cong
\widetilde{H}^{p}
\bigl(
(F,1_{\widehat {P(\Delta)}});K
\bigr)
$
for $F\in \Delta$ and all $p$,
it is a consequence of Theorem \ref{hochsterformula}.
\end{rem}

The following corollary of Theorem \ref{hochsterformula}
is an observation of
Stanley in \cite{ST87}.

\begin{cor}
\label{rational}
Let $\Sigma \subset \R^d$
be a rational pointed fan, let
$P$ be the face poset of $\Sigma$
and let $T$ be the $KP$-algebra considered
in Example \ref{toricrings}.
If $\Delta(P)$ is $K$-Cohen-Macaulay,
then $\lim T$ is Cohen-Macaulay.
\end{cor}
\begin{proof}
If $\Delta(P)$ is Cohen-Macaulay,
then
$$
\widetilde{H}^{p}\bigl((C,1_{\widehat{P}}); K\bigr)=0
\text{ for } C \in P \text{ and }
p \neq \dim \Delta(C,1_{\widehat{P}})
$$
where the intervals $(C,1_{\widehat{P}})$ are taken as subposets of $P$.
It is known that $P$ is a graded poset,
i.e. all maximal chains of $P$ have the same length.
We prove by induction on $\rank P-\rank C$ that
$
\dim \Delta(C,1_{\widehat{P}})=\dim (\lim T) - d_C-1
\text{ for }
C \in P$. This will conclude the proof by
the remark of \ref{hochsterformula}.

If $\rank C= \rank P$, then $C$ is a maximal face
and it is easy to see that
$\Ker (\lim T \to T_C)$ is a minimal prime ideal of $\lim T$.
Hence
$$
\dim \Delta(C,1_{\widehat{P}})
-1
=
\dim (\lim T) - \dim T_C-1
=
\dim (\lim T) - \dim C-1
$$
and we are done in this case.
If $\rank C<\rank P$, then
it is well known that there exists a face $C'\in P$
such that $\dim C'= \dim C+1$ and $C \subseteq C'$.
Thus by the induction hypothesis
$$
\dim \Delta(C,1_{\widehat{P}})
=
\dim \Delta(C',1_{\widehat{P}})+1
=
\dim (\lim T) - \dim C'
=
\dim (\lim T) - \dim C-1.
$$

\end{proof}

The next goal will be to find a result similar to \ref{cmring} for
Buchsbaum rings, at least in the $\Z$-graded case. Assume that
$K$ is a field and let $R=K[x_1,\dots,x_d]$ be an $\N$-graded
polynomial ring
(i.e. $R$ is $\Z$-graded such that $R_0=K$ and $R_i=0$ for $i<0$)
with graded maximal ideal $\mm=(x_1,\dots,x_d)$.
Recall that a finitely generated graded $R$-module $N$ is called a
{\em Buchsbaum module} if the $S_\mm$-module $N_\mm$ is Buchsbaum as
defined in \cite{STVO}.
If $K$ is an infinite field and $R$ is generated in degree $1$,
then $N$ is Buchsbaum if and only if $\dim N=0$ or $\dim N>0$ and every
homogeneous system of parameters $y_1,\dots,y_{\dim N}$ of $N$ is a
{\em weak $N$-sequence}, i.e. $(y_1,\dots,y_{i-1})N\colon y_i
=(y_1,\dots,y_{i-1}) N \colon \mm \text{ for } i=1,\dots,\dim N.$ If
$N$ is a Buchsbaum module, then $ \dim_K H_\mm^i(N)<\infty$ for
$i=1,\dots,\dim N-1$. In general the converse is not true. But if
there exists an $r \in \Z$ such that $H_\mm^i(N)_j=0$ for $j\neq r$
and $i=1,\dots,\dim N-1$, then $N$ is  Buchsbaum (see \cite[Satz
4.3.1]{SCH82}). For more details on the theory of Buchsbaum modules
we refer to Schenzel \cite{SCH82} and St\"uckrad and Vogel
\cite{STVO}.
\begin{cor}
\label{buchsring}
Let $R=K[x_1,\dots,x_d]$ be the $\N$-graded polynomial ring over a field $K$.
Let $P$ be a finite poset with a terminal element $0_P$ and
let $T$ be a cyclic flasque $RP$-algebra of $\Z$-graded Cohen-Macaulay $R$-algebras
such that $d_x=\dim T_x<d_y=\dim T_y$
for $x<y$ in $P$.
Assume that
$T_{0_P}=K$.
Then the following statements are equivalent:
\begin{enumerate}
\item
$\lim T$ is a Buchsbaum ring,
\item
for all
$x \in P \setminus \{0_P\}$
we have
$\widetilde{H}^{p}((x,1_{\widehat{P}});K)=0
\text{ for } p \neq \dim (\lim T) - d_x - 1.$
\end{enumerate}
\end{cor}
\begin{proof}
Note that $d_x>0$ for $x>0_P$.
If $\lim T$ is Buchsbaum,
then it follows from Theorem \ref{genhochster} that
for all
$x \in P$ such that $x \neq 0_P$
we have
$$
\Tilde{H}^{p}((x,1_{\widehat{P}});K)=0
\text{ for } p \neq \dim (\lim T) - d_x-1,
$$
because otherwise $H^i_\mm(\lim T)$
is not a finitely generated
$K$-vector space
since we have $\dim_K H^{d_x}_\mm(T_x)=\infty$ for $x \neq 0_p$.

If condition (ii) holds, then we have that
$H^i_\mm(\lim T)_j=0$ for $i<\dim (\lim T)$
and $j\neq 0$.
Thus $\lim T$ is a Buchsbaum ring.
\end{proof}

\section{$AP$-modules}
\label{kpmodules}
In the following sections we give
proofs of the results in Section \ref{mainresults}.
Let $A$ be an associative and unital ring and let
$P$ be a poset. A {\em left $AP$-module} $M$ is a system $(M_x)_{x
\in P}$ of left $A$-modules and homomorphisms $M_{xy} \colon M_y \to
M_x$ for $x \le y$ in $P$ with the property that $M_{xx}$ is the
identity on $M_x$ and that $M_{xy} \circ M_{yz} = M_{xz}$ for every
$x \le y \le z$ in $P$. A {\em homomorphism} $f \colon M \to N$ of
left $AP$-modules consists of homomorphisms $f_x \colon M_x \to N_x$
of left $A$-modules for $x \in P$ with the property that $f_x \circ
M_{xy} = N_{xy} \circ f_y$ for every $x \le y$ in $P$. We denote the
abelian group of homomorphisms from $M$ to $N$
by $\Hom_{AP}(M,N)$. The category-minded reader recognizes that a
left $AP$-module is a functor from $P$ to the category of left
$A$-modules, and that a homomorphism of left $AP$-modules is a
natural transformation of such functors. We denote the category of
left $AP$-modules by $AP\lmod$.

More generally, for every small category $\mathcal C$ enriched in
the category of abelian groups we could consider the category
$\mathcal C\lmod$ of enriched functors from $\mathcal C$ to the
category of abelian groups. Everything we do in this section and in
Section \ref{homologicalalg} in the category $AP\lmod$ can also be
done in the category $\mathcal C\lmod$. Since we are interested in
certain particular properties of the categories $AP\lmod$ we focus
on these.

A homomorphism $f \colon M \to N$ of left $AP$-modules is a {\em
monomorphism} if
$f_x \colon M_x \to N_x$ is injective for every $x \in P$, and it is
an {\em epimorphism}
if $f_x$
is surjective for every $x \in P$. A left $AP$-module $L$ is {\em
projective} if for every epimorphism $f \colon M \to N$ and every
homomorphism $g \colon L \to N$ of left $AP$-modules there exists a
homomorphism $\overline g \colon L \to M$ with $f \circ \overline g
= g$. Dually a left $AP$-module $I$ is {\em injective} if for every
monomorphism $f \colon M \to N$ and every homomorphism $g \colon M
\to I$ of left $AP$-modules there exists a homomorphism $\overline g
\colon N \to I$ with $\overline g \circ f = g$. The category
$AP\lmod$ of left $AP$-modules is an abelian category.
\begin{ex}
\label{projexample} \
\begin{enumerate}
  \item
  To an order-preserving map $f \colon P \to Q$ of posets and a left
  $AQ$-module $M$ there is an associated left $AP$-module $f^*M$ with
  $(f^*M)_x = M_{f(x)}$ and $(f^*M)_{xy} = M_{f(x)f(y)}$.
  We shall write $M$ instead
  of $f^*M$ when it is clear from the context that we are working with
  left $AP$-modules.
  \item
  If $Q$ is the one element poset, then the category of left $AQ$-modules
  is isomorphic to the category of left $A$-modules.
  \item
  Given a poset $P$ there is a unique order-preserving map $f \colon P \to Q$
  from $P$ to the one element poset $Q$. Thus we can consider a left
  $A$-module $E$ firstly as an $AQ$-module and secondly as an
  $AP$-module $f^*E$. Again, when it is clear from the context we
  write $E$ instead of $f^*E$. Note that $f^*E$ is the constant
  left $AP$-module
  with
  constant value $E$, that is, $f^* E_x = E$ for $x \in P$ and
  $f^*E_{xy}
  = \id_E$ for $x \le y$ in $P$.
\item
For $z \in P$ there is a {\em left $AP$-module $AP^z$
represented by $z$}. The left $AP$-module $AP^z$ takes $x$ to $AP^z_x
= A$ if $x \le z$ and to $0$ otherwise. The homomorphism
$AP^z_{xy}$ is the identity on $A$ if $x \le y \le z$ and
otherwise it is the zero homomorphism. If $M$ is another left
$AP$-module, then the abelian group of
homomorphisms from $AP^z$ to $M$ is isomorphic to the underlying
abelian group of $M_z$. In particular, $AP^z$ is a
projective left $AP$-module.
\item To a family $(M_i)_{i \in I}$ of left $AP$-modules we can
associate left $AP$-modules $\bigoplus_{i \in I} M_i$ and
$\prod_{i \in I} M_i$ with  $(\bigoplus_{i \in I} M_i)_x =
\bigoplus_{i \in I} (M_i)_x$ and
$(\prod_{i \in I} M_i)_x = \prod_{i\in I} (M_i)_x$ for $x \in P$.
\end{enumerate}
\end{ex}
\begin{defn}
The {\em limit} of a left $AP$-module $M$ is the left
$A$-submodule $\lim M$ of $\prod_{x \in P} M_x$ consisting of
sequences $m=(m_x)_{x \in P}$ with the property that $M_{xy}(m_y) =
m_x$ for every $x \le y$ in $P$.
\end{defn}
We call a left $AP$-module $M$
{\em finitely generated} is there exists an epimorphism of the form
$\bigoplus_{i \in I} AP^{z_i} \to M$ for some finite index set $I$.
Note that a finitely generated left $AP$-module is projective if and
only if it is a direct summand of
a finitely generated left $AP$-module
$\bigoplus_{i \in I} AP^{z_i}$ for some finite set $I$.
A left $AP$-module $M$ is called {\em Noetherian} if
every increasing sequence $M_1 \subseteq M_2 \subseteq \dots$ of
left sub-$AP$-modules of $M$ stabilizes. Recall that a poset $P$ is
topologically finite if $(0_{\widehat P},x]$ is
  finite for every $x \in P$.
\begin{lem}
\label{helper1}
  Suppose that $A$ is a left Noetherian ring and that $P$ is a topologically
  finite
  poset. Then $AP^z$
  is a Noetherian left $AP$-module for every $z \in P$.
\end{lem}
\begin{proof}
  A left sub-$AP$-module of $AP^z$ is uniquely determined by a family
  of left ideals $I(x)$ of $A$ for $x \le z$ in $P$ with the property
  that $I(y) \subseteq I(x)$ if $x \le y \le z$. By our assumptions
  this is a finite family of left ideals of a left Noetherian ring.
\end{proof}
It is well-known that if $A$ is left Noetherian, then every left
submodule of a finitely generated free $A$-module
is finitely
generated.
This also holds in the category of left $AP$-modules
in the following way:
\begin{prop}
\label{subfin}
If $AP^z$ is a Noetherian left $AP$-module for every $z \in P$, then
every submodule of a left $AP$-module
$\bigoplus_{i \in I} AP^{z_i}$  for some finite set $I$
is finitely generated.
\end{prop}
\begin{rem}
\label{incidencealg}
  If $P$ is finite and $A$ is a commutative ring, then we can
  consider the {\em incidence algebra}
$I(P,A)$
  of $P$ over $A$, that is, the
  $A$-algebra with underlying $A$-module
  $$I(P,A) =
  \bigoplus_{x \le y} A \cdot e_{x \le y}$$
  and with multiplication
  defined via $(e_{x \le y}) (e_{y \le z}) = e_{x \le z}$ and with $(e_{x
  \le y}) (e_{x' \le y'}) = 0$ if $y \ne x'$. (See
  \cite[Definition 3.6.1]{ST99} or \cite[p. 33]{Mitchell}.) Note
  that $\sum_{x \in P} e_{x
  \le x}$ is the multiplicative unit in $I(A,P)$ and that
  the elements
  $e_{x \le x}$ are idempotent in $I(A,P)$.

  Given a left module $M$ over
  the ring $I(A,P)$, we obtain $A$-modules $M_x := e_{x \le x}M$ for every $x
  \in P$, and $A$-linear homomorphisms $M_{xy} \colon M_y
  \to M_x$ given by
  $$(e_{y \le y}) m \mapsto (e_{x\le y})(e_{y \le y}) m = (e_{x\le
  y})m = (e_{x \le x}) (e_{x \le y})m$$
  for every
  $x \le y$ in $P$. In this situation
  the $A$-modules
  $(M_x)_{x \in P}$ and the $A$-homomorphisms $M_{xy}$ form a
  left $AP$-module. Conversely, if $M$ is a
  left $AP$-module, then the direct sum
  $\bigoplus_{x\in P} M_x$ of the $A$-modules can be given the
  structure of a left module over the ring $I(A,P)$ by defining $e_{x
  \le y} m_y = M_{xy} (m_y)$ for $m_y \in M_y$ and $e_{x \le y}m = 0$
  is $m \notin M_y$. This correspondence
  shows that the category of left modules over the ring $I(A,P)$ is
  equivalent to the category of left $AP$-modules. This
  justifies the above terminology and it shows that in the case where
  $P$ is finite the concept of $AP$-modules is really nothing
  new. However, as we shall see, many left $I(A,P)$-modules become
  more transparent when
  considered as left $AP$-modules.

Observe that the left $AP$-modules $AP^z$ correspond
to the ideals $I(P,A)e_{z\leq z}$.
The module $I(P,A)$ corresponds
to the left $AP$-module $\bigoplus_{z\in P} AP^z$.
\end{rem}

The next result is
well-known to specialists.
Since we did not find a proof in the literature we include it
for the sake of completeness.
\begin{prop}
\label{sheafpmod}
Let $P$ be a poset and let $A$ be an associative
and unital ring. The category of left $AP$-modules is equivalent to
the category of sheaves of left $A$-modules on $P$ with the
Alexandrov topology.
\end{prop}
\begin{proof}
Let $\Fc$ be a sheaf of left $A$-modules on $P$. We let $\Phi(P)$
denote the left $AP$-module with $\Phi(P)_x = \Fc((0_{\widehat
P},x])$ for $x \in P$ and with $\Phi(P)_{xy}$ equal to the
restriction homomorphism associated to the inclusion $(0_{\widehat
P},x]\subseteq (0_{\widehat P},y]$ for $x \le y$. This defines a
functor $\Phi$ from the category of sheaves of left $A$-modules on
$P$ to the category of left $AP$-modules.

Let $M$ be a left $AP$-module. For $U \subseteq P$ open define
\begin{equation}
\label{equation0}
\Psi(M)(U) = \{ (s_x) \in \prod_{x \in U} M_x : s_x=M_{xy}(s_y)
\text{ for } x\leq y\} = \lim M|_U.
\end{equation}
Given an inclusion $\iota \colon V \subseteq U$ of open subsets of $P$,
the restriction map $\Psi(M)(U) \to \Psi(M)(V)$ is the natural
restriction map of limits. Using that $(0_{\widehat P},x]$ is
contained in every neighborhood of
  $x$ it is straight forward to check that
$\Psi(M)$ is isomorphic to its associated sheaf of left $A$-modules
on $P$. (Compare  \cite[Proposition-Definition II.1.2]{RH}.)
Thus $\Psi(M)$ is a sheaf of left $A$-modules on $P$. We have
defined a functor $\Psi$ from the category of left $AP$-modules to
the category of sheaves of left $A$-modules on $P$.

If $M$ is a left $AP$-module, then the homomorphism
$$M_y \to \Phi(\Psi(M))_y = \lim M|_{(0_{\widehat P},y]}
$$
induced by the structure maps $M_{xy} \colon M_y \to M_x$ is an
isomorphism.

Conversely, if $\Fc$ is a sheaf of left $A$-modules on $P$ and $U$
is an open subset of $P$, then the intervals $(0_{\widehat P},x]$
for $x \in U$ form an open cover of $U$. Since the homomorphism
$$
\Fc((0_{\widehat P},y]) \to \Psi(\Phi(M))(0_{\widehat P},y]) =\lim_{x
\in (0_{\widehat P},y]} \Fc((0_{\widehat P},x])
$$
is an isomorphism for every $y \in P$, the sheaf condition on $\Fc$
ensures that the homomorphism
$$
\Fc(U) \to \Psi(\Phi(\Fc))(U)=\lim_{x \in U} \Fc((0_{\widehat P},x])
$$ is an isomorphism for every open subset $U$ of $P$.
This concludes the proof.
\end{proof}

We will also need to consider right $AP$-modules. A {\em right
$AP$-module} is a system $(M^x)_{x \in P}$ of right $A$-modules and
homomorphisms $M^{xy} \colon M^x \to M^y$ for $x \le y$ in $P$ with
the property that $M^{xx}$ is the identity on $M^x$ and that $M^{yz}
\circ M^{xy} = M^{xz}$ for every $x \le y \le z$ in $P$. A
homomorphism $f \colon M \to N$ of right $AP$-modules consists of
homomorphisms $f^x \colon M^x \to N^x$ of right $A$-modules for $x
\in P$ with the property that $f^y \circ M^{xy} = N^{xy} \circ f^x$
for every $x \le y$ in $P$. (In other words, a right $AP$-module is
a left $A^{\op} P^{\op} $-module.) The category $\rmod AP$ of right
$AP$-modules is also an abelian category.
\begin{ex}
For $x \in P$ there is a {\em right $AP$-module $AP_x$
represented by $x$}. The right $AP$-module $AP_x$ takes $z$ to $AP^z_x
= A$ if $x \le z$ and to $0$ otherwise. The homomorphism
$AP_x^{yz}$ is the identity on $A$ if $x \le y \le z$ and
otherwise it is the zero homomorphism. If $M$ is another left
$AP$-module, then the abelian group of
homomorphisms of right $AP$-modules from $M$ to $AP_x$ is isomorphic
to the underlying
abelian group of $M_x$. In particular, $AP_x$ is a
projective right $AP$-module.
\end{ex}

If $P = (P,\le)$ is a poset we let $(P,=)$ denote the poset with the
same elements as $P$ and with the partial order where no distinct
elements are comparable. Given a right $AP$-module $M$ and a left
$AP$-module $N$ the {\em tensor product $M \otimes_{A(P,=)} N$ of
$M$ and $N$ over
  $A(P,=)$} is
the abelian group $M \otimes_{A(P,=)} N = \bigoplus_{x \in P} M^x
\otimes_A N_x$. For every $x \in P$ there is a
homomorphism
$$N_* \colon AP_x
\otimes_{A(P,=)} N = \bigoplus_{y \in P} AP^y_x \otimes_{A} N_y \to
N_x$$ induced by the unique homomorphisms $AP^y_x \otimes_A N_y \to
N_x$ taking $1 \otimes n$ to $N_{xy}(n)$ for $n \in N_y$ and $x \le
y$. Similarly, there is a  homomorphism
$$M^* \colon M \otimes_{A(P,=)} AP^y
= \bigoplus_{x \in P} M^x \otimes_{A} AP^y_x \to M^y.$$ The {\em
tensor product $M\otimes_{AP} N$ of $M$ and $N$ over $AP$} is the
abelian group given by the cokernel of the homomorphism
\begin{displaymath}
  \bigoplus_{x,y \in P} M^x \otimes_A AP^y_x \otimes_A N_y
  \xto{M^* \otimes_A 1 - 1 \otimes_A N_*} \bigoplus_{x \in P} M^x
  \otimes_A N_x.
\end{displaymath}
Given posets $P$ and $Q$ and associative and unital rings $A$ and
$B$, an $AP$-$BQ$-bimodule is a system $(M_x^u)_{x \in P, u \in Q}$
of $A$-$B$-bimodules together with a left $AP$-module structure on
$(M_x^u)_{x \in P}$ for every $u \in Q$ and a right $BQ$-module
structures on $(M_x^u)_{u \in Q}$ for every $x \in P$ subject to the
condition that $M^v_{xy} \circ M^{uv}_y = M^{uv}_x \circ M_{xy}^u$
for $x \le y$ in $P$ and $u \le v$ in $Q$. A homomorphism $f \colon M \to N$ of
$AP$-$BQ$-bimodules consists of homomorphisms $f^u_x \colon M^u_x \to
N^u_x$ for $u \in Q$ and $x \in P$ such that for every $x \in P$ the
homomorphisms $(f^u_x)_{u \in Q}$ form a homomorphism of right
$BQ$-modules and for every $u \in Q$ the homomorphisms $(f^u_x)_{x \in
  P}$   form a homomorphism of left $AP$-modules. We denote by
$\Hom_{AP-BQ}(M,N)$ the abelian group of homomorphisms of
$AP$-$BQ$-bimodules
from $M$ to $N$.
If $P$ is a one-point poset,
then we say that $M$ is an $A$-$A'Q$-bimodule
instead of saying that it is an $AP$-$A'Q$-bimodule.
Similarly, if
$Q$ is a one-point poset,
then we say that $M$ is an $AP$-$A'$-bimodule.

\begin{ex}
\label{ex4}
For every poset $P$ and  every associative unital ring $A$ we can
consider the $AP$-$AP$-bimodule $AP$ with $AP_x^y = A$ if $x \le y$
and $AP_x^y = 0$ otherwise.
\end{ex}
Let $P,Q$ and $R$ be posets and let $A,A'$ and $A''$ be associative
and unital rings. If $M$ is an $A'Q$-$AP$-bimodule and $N$ is an
$AP$-$A''R$-bimodule, then the tensor product $M \otimes_{AP} N$
inherits the structure of an $A'Q$-$A''R$-bimodule.
Observe that if $P$ is a one-point poset, then
this is the $A'Q$-$A''R$-bimodule $M \otimes_{A} N$ induced by
$(M \otimes_{A} N)^{y}_{x}=M_x \otimes_{A} N^y$.

If further $L$ is an $AP$-$A'Q$-bimodule, then we let
$\Hom_{AP}(L,N)$ denote the $A'Q$-$A''R$-bimodule with
$\Hom_{AP}(L,N)_u^a$ given by the set of $AP$-homomorphisms from
$L^u = (L_x^u)_{x \in P}$ to $N^a = (N_x^a)_{x \in P}$, and with
structure homomorphisms induced from those of $L$ and $N$.
Observe that if $P$ is a one-point poset,
then this is the
$A'Q$-$A''R$-bimodule $\Hom_{A}(L,N)$ induced by
$\Hom_{A}(L,N)_u^a=\Hom_{A}(L^u,N^a)$.
Note also that
$\Hom_{AP}(L,N)$ is a left $A'$-module in the particular case where
$Q$ and $R$ are one-point posets and $A'' = \Z$.

Suppose that $A$ and $A'$ are algebras over a commutative ring $K$,
that is, there are given ring-homomorphism from $K$ to the centers
of $A$ and $A'$. If $M$ is a right $AP$-module and $M'$ is a right
$A'P'$-module, then the tensor product $M \otimes_K M'$ of $M$ and $M'$ over
$K$ is the right $(A \otimes_K A')(P \times P')$-module with $(M
\otimes_K M')^{(x,x')} = M^x \otimes_K M'^{x'}$. Similarly, if $N$
is a left $AP$-module and $N'$ is a left $A'P'$-module, then $N
\otimes_K N'$ is the left $(A \otimes_K A')(P \times P')$-module
with $(N \otimes_K N')_{(x,x')} = N_x \otimes_K N'_{x'}$.
\begin{prop}
\label{tensorhom}
  Suppose that $A$ and $A'$ are associative and unital algebras over a
  commutative ring $K$.
  If $M$ is a finitely generated projective left $AP$-module and $M'$ is a
  finitely generated projective left
  $A'P'$-module, then for every left $AP$-module $N$ and every left
  $A'P'$-module $N'$ the
  $\Hom$-$\otimes$ interchange homomorphism
  \begin{eqnarray*}
    \Hom_{AP}(M,N) \otimes_{K} \Hom_{A'P'}(M',N') &\to&
    \Hom_{(A\otimes_K A')(P
    \times P')}(M \otimes_{K} M',N \otimes_{K} N'), \\
    (f \otimes f') &\mapsto& (m \otimes m' \mapsto f(m) \otimes f'(m')).
  \end{eqnarray*}
  is an isomorphism.
\end{prop}
\begin{proof}
  Since every retract of an isomorphism is an isomorphism we can
  without loss of generality assume that $M$ and $M'$ are
  left $AP$-modules of the form $\bigoplus_{i \in I} AP^{z_i}$ for some finite set $I$.
  By direct sums, the statement is now reduced to the case $M = AP^y$
  and $M' =  A'P'^{y'}$.
In this case
$M \otimes_K M' \cong (A \otimes_K A')(P \times P')^{(y,y')}$.
We have that
\begin{eqnarray*}
&&\Hom_{AP}(M,N) \otimes_{K} \Hom_{A'P'}(M',N') \\
&=& \Hom_{AP}(AP^y,N) \otimes_{K} \Hom_{A'P'}(A'P'^{y'},N') \\
& \cong & N_y \otimes_{K} N'_{y'} = (N \otimes_K N')_{(y,y')} \\
&\cong&
\Hom_{(A\otimes_K A')(P\times P')}((A \otimes_K A')(P \times P')^{(y,y')},N \otimes_{K} N') \\
&\cong& \Hom_{(A\otimes_K A')(P
\times P')}(M \otimes_{K} M',N \otimes_{K} N'),
\end{eqnarray*}
where we used the fact (the Yoneda lemma)
that the homomorphisms
\begin{eqnarray*}
\Hom_{AP}(AP^y,N) &\to& N_y,\\
\varphi &\mapsto& \varphi_y(1_A)
\end{eqnarray*}
and
\begin{eqnarray*}
\Hom_{(A\otimes_K A')(P\times P')}((A \otimes_K A')(P \times P')^{(y,y')},N \otimes_{K} N') &\to& (N \otimes_K N')_{(y,y')},\\
\varphi &\mapsto& \varphi_{(y,y')}(1_{A \otimes_K A'})
\end{eqnarray*}
are isomorphisms.
\end{proof}
We leave the proof of the following lemma to the reader.
\begin{lem}
\label{hom-tens-adj}
  Let $A$, $A'$ and $A''$ be commutative rings and let $P$, $P'$ and
  $P''$ be posets. Given an $AP$-$A'P'$-bimodule $X$, an
  $A'P'$-$A''P''$-bimodule $Y$ and an $AP$-$A''P''$-bimodule $Z$ the
  homomorphism
  \begin{eqnarray*}
    \Hom_{AP-A''P''}(X \otimes_{A'P'} Y,Z) &\to&
    \Hom_{A'P'-A''P''}(Y,\Hom_{AP}(X,Z)) \\
    f &\mapsto& (y \mapsto (x \mapsto f(x\otimes y)))
  \end{eqnarray*}
  is an isomorphism.
\end{lem}

\section{Homological algebra}
\label{homologicalalg}
A {\em chain complex} $C$ of left
$AP$-modules is a collection $(C_n)_{n \in \Z}$ of left $AP$-module
together with homomorphisms $d = d_n \colon C_n \to C_{n-1}$ with
the property that $d \circ d = 0$. We call a chain complex $C$ {\em
positive} if $C_n = 0$ for $n<0$ and we call it {\em
  negative} if $C_n = 0$ for $n>0$.
The {\em homology} $H_*(C)$ of a chain complex $C$ is the collection
$(H_n(C))_{n \in \Z}$ of the left $AP$-modules $H_n(C) = \ker(d
\colon C_n \to C_{n-1})/\Image(d \colon C_{n+1} \to C_n)$. A {\em
homomorphism $f \colon C \to D$ of chain complexes of left
  $AP$-modules} is a collection of homomorphisms $f_n \colon C_n \to
D_n$ of left $AP$-modules with the property that $df_n = f_nd$. We
denote the category of chain complexes of left $AP$-modules
$\chcx(AP\lmod)$.

A {\em positive resolution} of a left $AP$-module $M$ is a positive
chain complex $C$ with the properties that $H_0(C) \cong M$ and that
$H_n(C) = 0$ for $n \ne 0$. A {\em projective resolution} of $M$ is
a positive resolution $C$ of $M$ with the property that $C_n$ is a
projective left $AP$-module for every $n$.
\begin{prop}
\label{helper2}
  If $AP^z$ is a Noetherian left $AP$-module for every $z \in P$, then
  there exists a (homological) degreewise finitely generated projective resolution
  of every finitely generated left $AP$-module $M$.
\end{prop}
\begin{proof}
  There exists a short exact sequence of
  the form
  \begin{displaymath}
    0 \to G_0 \to F_0 \to M \to 0,
  \end{displaymath}
  where $F_0$ is a finitely generated left $AP$-module of the form $\bigoplus_{i \in I} AP^{z_i}$ for some finite set $I$.
  By Proposition \ref{subfin} $G_0$ is finitely generated, and
  proceeding by induction we obtain a resolution $F \to M$, with $F_i$
  a finitely generated and projective left $AP$-module for every $i \ge 0$.
\end{proof}

A {\em negative resolution} of a left $AP$-module $M$ is a negative
chain complex $C$ with the properties that $H_0(C) \cong M$ and that
$H_n(C) = 0$ for $n \ne 0$. An {\em injective resolution} of $M$ is
a negative resolution of $M$ consisting of injective left
$AP$-modules. Every left $AP$-module has a projective resolution, so
the category $AP\lmod$ is an abelian category with enough
projectives. The usual argument
 showing that module categories are categories with enough
injectives (see  \cite{WE95}) also shows that every left
$AP$-module has an injective resolution.
This fact can also be deduced from
Proposition \ref{sheafpmod} and the fact that the category
of sheaves of left $A$-modules on $P$
has enough injective objects
(see \cite[Proposition 2.2]{RH}).
The above observations also
apply to right $AP$-modules and to $AP$-$LQ$-bimodules.

Suppose that $A$ is an algebra over a commutative ring $K$. The
functor
\begin{displaymath}
 \otimes_{AP} \colon AP\lmod \times \rmod AP \to K \lmod, \quad (M,N)
 \mapsto M \otimes_{AP} N
\end{displaymath}
is additive and right exact in both $M$ and $N$ and it can be
extended to a functor
\begin{displaymath}
 \otimes_{AP} \colon \chcx(AP\lmod) \times \chcx(\rmod AP) \to
 \chcx(K \lmod), \quad (C,D)
 \mapsto C \otimes_{AP} D.
\end{displaymath}
More precisely, $(C \otimes_{AP} D)_n = \bigoplus_{r+s = n} C_r
\otimes_{AP} D_s$ and if $m \in (C_r)^x$ and $n \in (D_s)_x$, then the
differential takes the element in $C_r \otimes_{AP} D_s$ represented
by $m \otimes n \in (C_r)^x \otimes_A (D_s)_x$ to the sum of the
elements in $(C \otimes_{AP} D)_{n-1}$ represented by $dm \otimes n$
and $(-1)^rm \otimes dn$.

If $C_M$ is a projective resolution of $M$ and $C_N$ is a projective
resolution of $N$, then the $\Tor^{AP}_n(M,N)$ denotes
the isomorphism class of
\begin{displaymath}
  H_n(C_M \otimes_{AP} N) \cong H_n(M \otimes_{AP} C_N).
\end{displaymath}
It is a well-known fact of homological algebra that the above
homology groups are universal $\delta$-functors and that they are
naturally isomorphic by a unique isomorphism. See for example Weibel
\cite{WE95}. Similarly, the functor
\begin{displaymath}
  \Hom_{AP} \colon AP\lmod \times AP\lmod \to K\lmod, \quad (M,N)
  \mapsto \Hom_{AP}(M,N)
\end{displaymath}
is additive in both its entries.
It is left exact in $M$ and in $N$ and it can be extended to a functor
\begin{displaymath}
  \Hom_{AP} \colon \chcx(AP\lmod) \times \chcx(AP\lmod) \to
  \chcx(K\lmod), \quad (C,D)
  \mapsto \Hom_{AP}(C,D)
\end{displaymath}
with $\Hom_{AP}(C,D)_n = \prod_{-r+s = n} \Hom_{AP}(C_r,D_s)$ and
with $df$ given by $(df)(c) = f(dc) - d(f(c))$ for $f \in
\Hom_{AP}(C,D)_n$. If $D_M$ is a projective resolution of $M$ and
$C_N$ is an injective resolution of $N$, then the functors
$\Ext_{AP}^n(M,N)$ denotes the isomorphism class of
\begin{displaymath}
  \Ext_{AP}^n(M,N) = H_{-n}(\Hom_{AP}(D_M, N)) = H_{-n}(\Hom_{AP}(M,C_N)).
\end{displaymath}
Again it is a well-known fact of homological algebra that the above
homology groups are universal $\delta$-functors and that they are
naturally isomorphic by a unique isomorphism.

\begin{rem}
\label{remoncoh}
Given a sheaf $\Fc$ of left $A$-modules let
$H^n(P,\Fc)$ denote the sheaf cohomology. For a left $AP$-module $M$
we have that
$$
\Ext^0_{AP}(A,M)\cong \Hom^0_{AP}(A,M) \cong H^0(P,\Psi(M))
$$
where $\Psi(M)$ is the sheaf associated to $M$ as constructed in (\ref{equation0})
in the proof of \ref{sheafpmod}.
Using Proposition \ref{sheafpmod} we see that if $M$ is injective
then $\Psi(M)$ is an injective sheaf on $P$. It follows that there
is a natural isomorphism
$$\Ext^n_{AP}(A,M) \cong H^n(P,\Psi(M))$$
for every $n$. The above groups are isomorphic to the
Hochschild--Mitchell cohomology groups $H^n(P,M)$ of $P$ with
coefficients in $M$. More precisely, the Hochschild--Mitchell
complex is the chain complex $\Hom_{AP}(B^K(AP,AP,A),M)$, where
$B^K(AP,AP,A)$ is a particular projective resolution of $A$ over
$AP$ called the bar-construction
(see \cite{BW} ,\cite{Mitchell}).
If $M$ is a left $A$-module considered as a constant $AP$-module, that is,
$M_{xy}$ is the identity on $M$ for every $x \le y$ in $P$, then we
can consider the cohomology $H^n(\Delta(P),M)$
of the simplicial complex $\Delta(P)$ with coefficients in $M$. In
this case the chain complex $\Hom_{AP}(B^K(AP,AP,A),M)$ is
isomorphic to the chain complex computing $H^n(\Delta(P),M)$. Thus
there are natural isomorphisms
\begin{equation}
\label{equation}
H^n(\Delta(P),M) \cong H^n(P,M) \cong \Ext^n_{AP}(A,M) \cong
H^n(P,\Psi(M)).
\end{equation}
Furthermore McCord has shown that the cohomology
groups $H^n(P,M)$ are isomorphic to the singular cohomology groups
of $P$ with the Alexandrov topology and with coefficients in $M$
\cite{Mccord}.
\end{rem}
Borrowing notation from the theory of sheaves
we call a left $AP$-module $M$ {\em flasque}
if $\lim M|_U \to \lim M|_V$
is surjective for all open sets
$V \subseteq U$ of the poset $P$.
Note that
$\Hom_{AP}(A,M) \cong \lim M$ is isomorphic to the group of global
sections  of the sheaf associated to $M$ under the equivalence of
categories between left $AP$-modules and sheaves of $A$-modules on $P$ of
Proposition \ref{sheafpmod}.
Since the higher cohomology groups of a
flasque sheaf vanish (see \cite[Proposition III,2,5]{RH}) we have:
\begin{lem}
\label{flasquecrit0}
Let $A$ be an associative and unital ring and let $P$ be a poset.
If $M$ is a flasque left $AP$-module, then
$\Ext_{AP}^i(A,M)=0 \text{ for all } i>0.$
\end{lem}

\begin{prop}
\label{hom-tensor}
  Suppose that $A$ and $A'$ are associative and unital algebras over a
  commutative ring $K$.
  If $C$ is a positive chain complex of finitely generated projective left
  $AP$-modules and $C'$ is a positive chain complex of finitely
  generated projective left
  $A'P'$-modules, then for every negative chain complex $D$ of left
  $AP$-modules
  and every negative chain complex $D'$ of left
  $A'P'$-modules the
  $\Hom$-$\otimes$ interchange homomorphism
  \begin{eqnarray*}
    \Hom_{AP}(C,D) \otimes_{K} \Hom_{A'P'}(C',D') &\to&
    \Hom_{(A\otimes_K A')(P \times P')}(C \otimes_{K} C',D \otimes_{K} D'), \\
    (f \otimes f') &\mapsto& (m \otimes m' \mapsto (-1)^{|f'||m|} f(m)
    \otimes f'(m')).
  \end{eqnarray*}
  is an isomorphism of chain complexes.
\end{prop}
\begin{proof}
  This is a direct consequence of Proposition \ref{tensorhom}.
\end{proof}
Note that
Lemma \ref{hom-tens-adj} also holds for chain complexes.
\section{Local Cohomology}
\label{locorp}
Let $R$ be a commutative ring, let $I$ be a finitely
generated ideal in $R$ and let $N$ be an $R$-module. The natural
projections $R/I^{n+1} \to R/I^{n}$ induce  maps
$\Ext^q_R(R/I^{n},N) \to \Ext^q_R(R/I^{n+1},N)$.
Our model for the
{\em $q$-th local cohomology group of $N$ with respect to $I$}
is the colimit
$$
H^q_I(N) = \colim \Ext^q_R(R/I^n,N).
$$

\begin{prop}
\label{prop1}
Let $I$ be a finitely generated ideal of a commutative ring $R$,
let $P$ be a poset and let $M$ be a left $RP$-module.
Suppose that there exist
a degreewise finitely generated projective
resolution  of $R$ over $RP$
and
a degreewise finitely generated
free resolution of $R/I^n$ over $R$.
If
$\Ext^q_{RP}(R,M) = 0$ for $q > 0$,
then there is a
natural isomorphism
$$
H^q_I(\lim M) \cong \colim \Ext^q_{RP}(R/I^n,M)
\text{ of $R$-modules for every } q \ge 0.
$$
\end{prop}
\begin{proof}
Let $E$ be a degreewise finitely generated projective
resolution of $R$ over $RP$
and let
$F_n \to R/I^n$ be a degreewise finitely generated
free resolution  of $R/I^n$ over $R$.
We may consider $F_n \otimes_R E$ as a projective resolution of
$R/I^n$ over $RP$. The vanishing of $\Ext^q_{RP}(R,M)$ implies that
the homomorphism
$E \to R$ induces a quasi-isomorphism $\Hom_{RP}(R,M) \to
\Hom_{RP}(E,M)$. Applying this quasi-isomorphism,
$\otimes$-$\Hom$-interchange and basic isomorphisms of the form
$\lim M \cong \Hom_R(R,\lim M)$ and $R \otimes_R N \cong N$ and
noting that $\Hom_R(F_n,R)$ is a degreewise free $R$-module, we
obtain the following chain of isomorphisms and quasi-isomor\-phisms:
  \begin{eqnarray*}
    \Hom_{RP} (F_n \otimes_R E,M)
    &\cong&
    \Hom_{R \otimes_R R(* \times P)} (F_n \otimes_R E, R \otimes_R M) \\
    &\cong&
    \Hom_R(F_n,R) \otimes_R \Hom_{RP}(E,M) \\
    &\simeq&
    \Hom_R(F_n,R) \otimes_R  \Hom_{RP}(R,M)\\
    &\cong&
    \Hom_R(F_n,R) \otimes_R \lim M \\
    &\cong&
    \Hom_R(F_n,R) \otimes_R \Hom_R(R,\lim M) \\
    &\cong&
    \Hom_R(F_n \otimes_R R,R \otimes_R \lim M ) \\
    &\cong&
    \Hom_R(F_n,\lim M).
  \end{eqnarray*}
Taking cohomology we get the natural isomorphism
\begin{displaymath}
  \Ext^q_{RP}(R/I^n,M) \cong \Ext^q_{R}(R/I^n,\lim M)
\end{displaymath}
of $R$-modules. Forming the colimit of these isomorphism we obtain
the isomorphism
$$
\colim \Ext^q_{RP}(R/I^n,M) \cong \colim
\Ext^q_{R}(R/I^n,\lim M) = H^q_I(\lim M).
$$
\end{proof}
In our applications we need a graded version of
the above result. If $R$ is a $\Z^d$-graded commutative ring and
$N$ and $N'$
are $\Z^d$-graded $R$-modules, then the group
$\Hom_R^{\mathrm{gr}}(N,N')$ of
homogeneous
homomorphisms from $N$ to $N'$ is a $\Z^d$-graded $R$-module.
Choosing a
projective resolution $E$ of $N$ in the category of $\Z^d$-graded
$R$-modules and
homogeneous homomorphisms of degree zero we obtain a chain complex
$\Hom_R^{\mathrm{gr}}(E,N')$ of $\Z^d$-graded $R$-modules.

If $F$ is a finitely generated free $\Z^d$-graded $R$-module, then
$\Hom_R(F,N')$ is isomorphic to $\Hom_R^{\mathrm{gr}}(F,N')$ for every
$\Z^d$-graded $R$-module $N'$. Since both $\Hom_R(-,N')$ and
$\Hom_R^{\gr}(-,N')$ are left exact functors it follows for every
finitely presented $\Z^d$-graded $R$-module $N$ that
$\Hom_R^{\mathrm{gr}}(N,N')$ and $\Hom_R(N,N')$ are isomorphic for every
$\Z^d$-graded $R$-module $N'$.
In particular, a finitely generated projective $\Z^d$-graded
$R$-module is also  projective considered as a non-graded $R$-module.
If $R$ is Noetherian and $N$ is a finitely generated $\Z^d$-graded $R$-module,
we
obtain a $\Z^d$-grading of
\begin{displaymath}
  \Ext^{q}_R(N,N') \cong H_{-q}(\Hom(E,N'))
\end{displaymath}
for every $q \ge 0$.
Here $E$ is a degreewise finitely generated projective resolution of
the $\Z^d$-graded $R$-module $N$.
In this case we obtain a grading of $H^q_I(N)$
if $I$ is a finitely generated graded ideal in $R$.

Recall that
for a poset $P$  a $\Z^d$-graded left $RP$-module is a left $RP$-module $M$ together
with gradings of the $R$-modules $M_x$ such that the homomorphisms
$M_{xy} \colon M_y \to M_x$ are homogeneous of degree zero for every
$x \le y$ in $P$.
The proof of Proposition \ref{prop1} can easily be modified to a proof
of the following result:
\begin{prop}
\label{prop12}
Let $I$ be a finitely generated graded ideal of a commutative $\Z^d$-graded ring $R$,
let $P$ be a poset and let $M$ be a $\Z^d$-graded left $RP$-module.
Suppose that there exist a degreewise finitely generated $\Z^d$-graded projective
resolution of $R$ over $RP$
and
a degreewise finitely generated $\Z^d$-graded
free resolution of $R/I^n$ over $R$.
If $\Ext^q_{RP}(R,M) = 0$ for $q > 0$,
then there is a
natural isomorphism
$$
H^q_I(\lim M) \cong \colim \Ext^q_{RP}(R/I^n,M)
$$
of $\Z^d$-graded $R$-modules for every $q \ge 0$.
\end{prop}

\begin{defn}
  Let $P$ be a poset and $R$ a commutative ring. For $x \in P$ and $q
  \in \Z$, the {\em left $RP$-skyscraper chain complex $R(x,q)$} is the left
  $RP$-module with $R(x,q)_y = 0$ for $y \ne x$ and with $R(x,q)_x =
  R[q]$ equal to the chain complex consisting of a copy of $R$ in
  homological degree $-q$.
\end{defn}

In the following
a zig-zag
chain of quasi-isomorphisms between complexes
$C$ and $D$
consists of
chain complexes $E_0,\dots,E_{2n}$ with $E_0=C$, $E_{2n}=D$,
of quasi-isomorphisms $E_{2k-1} \to E_{2k-2}$
and
of quasi-isomorphisms $E_{2k-1} \to E_{2k}$
for $k=1,\dots,n$.

Recall the discussion about bimodules after Example \ref{ex4}.
For example we need the following.
Assume that $K \to R$ is a homomorphism of commutative rings,
$F$ an $R$-module,
$P$ a poset and $M$ a left $RP$-module.
Then $M$ is also an $R-KP^{\op}$-bimodule. Hence
$\Hom_R(F,M)$ has a $R-KP^{\op}$-bimodule structure
and thus a left $RP$-module structure
induced by $\Hom_R(N,M)_x=\Hom_R(N,M_x)$ for $x \in P$.
This gives also an $RP$-module structure on $\Ext_R(N,M)$.

\begin{thm}
\label{decomp}
  Let $K \to R$ be a homomorphism of commutative rings,
  let $I$ be a finitely generated ideal of $R$,
  let $P$ be a poset and let $M$ be a left $RP$-module.
  Suppose that $E$ is a positive chain complex of
  finitely generated projective left $KP$-modules
  and that we have
  degreewise finitely generated
  free resolutions $F_n \to R/I^n$ of $R/I^n$ over $R$
  together with homomorphisms $F_{n+1} \to F_n$
  inducing the natural
  projections $R/I^{n+1} \to R/I^n$ in homology.
  Assume that the following
  conditions are satisfied:
  \begin{enumerate}
    \item
    there exists a zig-zag chain of quasi-isomorphisms between the
    chain complex of left $KP$-modules
    $\colim \Hom_R(F_n,M)$
and the chain complex of left $KP$-modules $H_I^{-*}(M):= \colim \Ext^{-*}_{R}(R/I^n,M)$ with
    $H_I^{-*}(M)_x = H_I^{-*}(M_x)$ for $x \in P$,
    \item
    for every $x < y$ in $P$ the homomorphism
    $H_I^{-*}(M_y) \to H_I^{-*}(M_x)$
    is the
    zero-homo\-mor\-phism.
  \end{enumerate}
  Then there is a
  natural zig-zag chain of quasi-isomorphisms of chain complexes of
  $K$-modules of the form
  \begin{displaymath}
    \colim \Hom_{RP}(F_n \otimes_K E,M
    ) \simeq
    \bigoplus_{x \in P}
    \ \ \bigoplus_{q \ge 0}  \Hom_{KP}(E,K(x,q)) \otimes_K
    H_I^q(M_x).
  \end{displaymath}
\end{thm}
\begin{proof}
  The asserted weak equivalence is the composition of the following
  homomorphisms:
\begin{eqnarray*}
   \colim \Hom_{RP}(F_n \otimes_K E, M) &\cong& \colim
   \Hom_{R-KP^{\op}}(F_n \otimes_K E,M) \\
   &\cong& \colim \Hom_{K - KP^{\op}}(E,\Hom_R(F_n,M)) \\
   &\cong& \Hom_{KP}(E,\colim \Hom_R(F_n,M)) \\
   &\simeq& \Hom_{KP}(E,H^{-*}_I(M)) \\
   &\cong& \Hom_{KP}(E,\bigoplus_{x \in P} \, \bigoplus_{q \ge 0}
   K(x,q) \otimes_K H_I^q(M)) \\
   &\cong& \bigoplus_{x \in P} \, \bigoplus_{q \ge 0}\Hom_{KP}(E,
   K(x,q) \otimes_K H_I^q(M)) \\
   &\cong& \bigoplus_{x \in P} \, \bigoplus_{q \ge 0}\Hom_{KP}(E,
   K(x,q)) \otimes_K H_I^q(M).
 \end{eqnarray*}
  Here the first isomorphism is clear,
  the second isomorphism is given by Lemma
  \ref{hom-tens-adj} and the third isomorphism is due to the facts that
  $E$ is positive and degreewise finitely generated and that
  $\Hom_R(F_n,M)$ is negative. The zig-zag of quasi-isomorphism on the fourth
  line exists by condition (i) and the fact that $\Hom_{KP}(E,-)$
  preserves quasi-isomorphisms between negative chain complexes. The
  isomorphism on the fifth line is a direct consequence of condition
  (ii).
  The isomorphism on the sixth line is again due to the fact that $E$
  is positive
  and degreewise finitely generated and that $K(x,q)$ is concentrated
  in one homological degree. The last isomorphism is a direct
  application of Proposition \ref{hom-tensor}.
\end{proof}
For reference we state a $\Z^d$-graded version of the above result. It is
proved in exactly the same way.
\begin{thm}
\label{decomp2}
  Let $K \to R$ be a homomorphism of $\Z^d$-graded commutative
  rings. Let $I$ be a finitely generated graded ideal in $R$,
  let $P$ be a poset and let $M$ be a $\Z^d$-graded left $RP$-module.
  Suppose that $E$ is a positive chain complex of
  finitely generated $\Z^d$-graded projective left $KP$-modules
  and that we have
    degreewise finitely generated $\Z^d$-graded
    free resolutions $F_n \to R/I^n$ of $R/I^n$ over $R$
  together with homomorphisms $F_{n+1} \to F_n$
  inducing the natural
  projections $R/I^{n+1} \to R/I^n$ in homology.
  Assume that the following
  conditions are satisfied:
  \begin{enumerate}
    \item
    there exists a zig-zag chain of homogeneous quasi-isomorphisms of
    degree zero between the chain complex of
    left $KP$-modules
    $\colim \Hom_R(F_n,M)$
    and the chain complex of left $KP$-modules $H_I^{-*}(M)$ with
    $H_I^{-*}(M)_x =
    H_I^{-*}(M_x)$ for $x \in P$,
    \item
    for every $x < y$ in $P$ the homomorphism
    $H_I^{-*}(M_y) \to H_I^{-*}(M_x)$
    is the
    zero-homo\-mor\-phism.
  \end{enumerate}
  Then there is a
  natural zig-zag chain of homogeneous quasi-isomorphisms of chain
  complexes of $\Z^d$-graded $K$-modules of the form
  \begin{displaymath}
    \colim \Hom_{RP}(F_n \otimes_R E,M) \simeq \bigoplus_{x \in P}
    \ \ \bigoplus_{q \ge 0}  \Hom_{KP}(E,K(x,q)) \otimes_K
    H_I^q(M_x).
  \end{displaymath}
\end{thm}

\begin{prop}
\label{helper3}
  Let $P$ be a finite poset, let $K$ be a field and let $C$ be a chain
  complex of left $KP$-modules.
  Suppose that for every $x\in P$
  there exists $n_x \in \Z$ such that $H_*(C_x)$ is concentrated in
  degree $n_x$ and that $x < y$ implies $n_x > n_y$. Then there exists
  a zig-zag chain of quasi-isomorphisms between
  the chain complexes of left $KP$-modules
  $C$ and $H_*(C)$.
\end{prop}
\begin{proof}
Let $\partial$ be the differential of $C$ and
$$B_k(C)=\Image(\partial_{k+1}) \subseteq C_k,\quad
Z_k(C)=\Ker(\partial_{k}) \subseteq C_k$$
for all $k$.
We define a chain complex $C'$ of left $KP$-modules with
  \begin{displaymath}
    (C'_x)_k =
    \begin{cases}
      (C_x)_k & k > n_x, \\
      B_k(C_x) \oplus H_k(C_x) & k = n_x, \\
      0 & k < n_x.
    \end{cases}
  \end{displaymath}
  For $k > n_x + 1$ the boundary $d \colon (C'_x)_{k} \to
  (C'_x)_{k-1}$ is the boundary map of $C_x$, for $k = n_x + 1$ it is
  the map $(d,0) \colon (C_x)_{n_x + 1} \to B_{n_x}(C_x) \oplus
  H_{n_x}(C_x)$, and for $k \le n_x$ it is the zero map. We give $C'$
  the structure of a chain complex of left $KP$-modules, where the map
  $(C'_{xy})_k \colon (C'_y)_k \to (C'_x)_k$
  is the map
  $(C_{xy})_k \colon (C_y)_k \to (C_x)_k$
  for $k > n_x$, the zero map for $k < n_x$ and
  the map
  $(C_y)_k \xto {(C_{xy})_k} (C_x)_k \to B_{k}(C_x) \oplus H_{k} (C_x)$
  for $k = n_x$.
  Here
  $(C_x)_{n_x} \to B_{n_x}(C_x)$ is a (chosen) retract
  of the inclusion $B_{n_x}(C_x) \to (C_x)_{n_x}$.
  Similarly, the map
  $(C_x)_{n_x} \to H_{n_x}(C_x)$ is the composition $(C_x)_{n_x} \to
  Z_{n_x}(C_x) \to H_{n_x}(C_x)$, where $(C_x)_{n_x} \to
  Z_{n_x}(C_x)$ is
  a (chosen) retract
  of the inclusion $Z_{n_x}(C_x) \to (C_x)_{n_x}$.
  Note that these retracts exist since
  $K$ is a field.

  There is a quasi-isomorphism $f \colon C \to C'$ with $(f_x)_k$
  equal to the identity on $(C_x)_k$ if $k > n_x$, the zero map if $k
  < n_x$ and the map $(C_x)_k \to B_{k}(C_x) \oplus H_{k} (C_x)$
  defined above for
  $k = n_x$.

  On the other hand, the inclusion $H_*(C) \to C'$ is a
  quasi-isomorphism of chain
  complexes of left $KP$-modules. Thus there are
  quasi-isomorphisms $C \to C' \gets H_*(C)$.
\end{proof}
\begin{lem}
\label{helper4}
  Let $P$ be a poset, let $x \in P$ and let $q \ge 0$. For every $n
  $ there is a
  natural isomorphism
  \begin{displaymath}
    \widetilde H^{n-q-1}((x,1_{\widehat P});R) \cong \Ext^n_{RP}(R,R(x,q)).
  \end{displaymath}
\end{lem}
\begin{proof}
  Let $F = RP^y$ for $y \in P$, that is, $F_x = R$ if $x \le y$ and
  $F_x = 0$ otherwise. Given a left $RP$-module $M$,
  the so-called Yoneda lemma
  provides an isomorphism
  \begin{displaymath}
   \Hom_{R(x,1_{\widehat P})}(F|_{R(x,1_{\widehat
  P})},M|_{R(x,1_{\widehat P})}) \cong
   \begin{cases}
    M_y & \text{ if } x < y, \\
    0 & \text{ otherwise.}
   \end{cases}
  \end{displaymath}
  If $x < y$ then the above isomorphism takes $\varphi \colon
  F|_{R(x,1_{\widehat
  P})} \to M|_{R(x,1_{\widehat P})}$ to $\varphi_y(1)$, where $1$ is
  the unit of $R = RP^y_y$.
  Similarly, there are isomorphisms
  \begin{displaymath}
   \Hom_{R[x,1_{\widehat P})}(F|_{R[x,1_{\widehat
  P})},M|_{R[x,1_{\widehat P})}) \cong
   \begin{cases}
    M_y & \text{ if } x \le y, \\
    0 & \text{ otherwise, }
   \end{cases}
  \end{displaymath}
  and
  \begin{displaymath}
   \Hom_{RP}(F,R(x,q)) \cong
   \begin{cases}
    R[q] & \text{ if } x = y, \\
    0 & \text{ if }  y \ne x.
   \end{cases}
  \end{displaymath}
  Let $E$ be a projective resolution of $R$ as a left $RP$-module.
  Since the above isomorphisms are natural in the projective left
  $RP$-module $F$,
  there is a short exact sequence of left $RP$-modules of the form
  \begin{displaymath}
    0 \to \Hom_{R(x,1_{\widehat P})}(E|_{R(x,1_{\widehat P})},R[q]) \to
    \Hom_{R[x,1_{\widehat P})}(E|_{R[x,1_{\widehat P})},R[q]) \to
    \Hom_{RP}(E,R(x,q)) \to 0.
  \end{displaymath}
  Note that $E|_{R(x,1_{\widehat P})}$ is a projective resolution of $R$
  over $R(x,1_{\widehat P})$ and
  that $E|_{R[x,1_{\widehat P})}$ is a projective resolution of $R$
  over $R[x,1_{\widehat P})$. Thus the long exact sequence associated
  to the above short exact sequence of chain complexes is of the form
  \begin{displaymath}\\
    \Ext^{n-1}_{RP}(R,R(x,q)) \to \Ext^{n-q}_{R(x,1_{\widehat
    P})}(R,R) \to
    \Ext^{n-q}_{R[x,1_{\widehat P})}(R,R) \to
    \Ext^n_{RP}(R,R(x,q)) \to
    \dots
  \end{displaymath}
  The result now follows from the fact that $[x,1_{\widehat P})$ is
 contractible
 since by Equation (\ref{equation}) in Remark \ref{remoncoh} we have $\Ext^{j}_{RQ}(R,R) \cong
 H^j(Q,R)$ for any poset $Q$.
\end{proof}

Now we are able to present the proof of
Theorem \ref{genhochster02}.

\begin{proof}[Proof of Theorem \ref{genhochster02}]
Since $P$ is finite,
$K$ is a finitely generated $KP$-module
and
it follows from Lemma \ref{helper1} and Proposition \ref{helper2}
that $K$ has a
degreewise finitely generated projective
resolution $E$ over $KP$.
All the rings $R/I^k$ have free
resolution $F_k$ over $R$ which are
degreewise finitely generated.
By Proposition \ref{prop1}
there is an isomorphism
$$
H^i_I(\lim T)
\cong
\colim \Ext^i_{RP}(R/I^k,T).
$$
By the assumption (i)
on $T_x$ we have
$H^i_I(T_x)=0$ for $i \neq d_x$.
By Proposition \ref{helper3}
there exists a zig-zag chain of  quasi-isomorphisms
between the left
$KP$-modules
$\colim \Hom_R(F_k,T)$
and
$$H^{\pnt}_I(T) := \colim \Ext^{-\pnt}_{R}(R/I^k,T) \cong
H_{\pnt}(\colim \Hom_R(F_k,T)).$$
(Here we used the fact that homology and filtered colimits commute.)
Since $x < y$ implies
$d_x<d_y$
the ho\-momorphism $H^i_I(T_y) \to H^i_I(T_x)$
is the zero-ho\-mo\-morphism.
Hence Theorem \ref{decomp}
implies that there is a
natural zig-zag chain of
quasi-isomorphisms of the form
  \begin{displaymath}
    \colim \Hom_{RP}(F_n \otimes_K E,T
    ) \simeq
    \bigoplus_{x \in P}
    \ \ \bigoplus_{q \ge 0}  \Hom_{KP}(E,K(x,q)) \otimes_K
    H_I^q(T_x).
  \end{displaymath}

Since $K$ is a field we obtain the isomorphism
$$
H^i_I( \lim T) \cong \colim
\Ext_{RP}^{i}(R/I^k,T) \cong
\bigoplus_{x \in P}
\ \ \bigoplus_{0 \le q \le i}
\Ext_{KP}^{i}(K,K(x,q)) \tensor_K
H^q_I(T_x)
$$
of  $K$-modules.
Lemma \ref{helper4}
implies that
$$
H^i_I( \lim T)
\cong
\bigoplus_{x \in P}
\ \ \bigoplus_{0 \le q \le i}
\widetilde H^{i-q-1}((x,1_{\widehat P});K)
\tensor_K
H^q_I(T_x)
$$
By the assumption (i) we obtain
\begin{eqnarray*}
H^i_I(\lim T)
&\cong&
\bigoplus_{x \in P}
\widetilde H^{i - d_x -1}((x,1_{\widehat P});K)
\tensor_K
H^{d_x}_I(T_x).
\end{eqnarray*}
The $\Z^d$-graded version can be proved in exactly the same way
since Proposition \ref{prop1} and Theorem \ref{decomp}
have $\Z^d$-graded versions.
This concludes the proof.
\end{proof}

\begin{rem}
Some of the results in this section can be interpreted in terms of Grothen\-dieck
 spectral sequences. Given abelian categories $\Ac$, $\Bc$ and $\Cc$
 with enough injectives and left exact functors $G \colon \Ac \to
 \Bc$ and $F \colon \Bc \to \Cc$ the Grothendieck spectral sequence
 with $E_2$-term $E_2^{p,q}(A) = (R^pF)(R^qG(A))$ given by composed
 right derived functors converges to the right derived functor
 $(R^{p+q}(FG))(A)$ of $FG$ for every object $A$ of $\Ac$.

 The $q$-th local cohomology $H^q_I(N)$ of an $R$-module
 $N$ is the $q$-th right derived of the zeroth local cohomology
 $H^0_I(N) \colon R\lmod \to R\lmod$ and the group $\Ext^q_{RP}(R,M)$
 is the $q$-th right
 derived functor of the functor $\lim_{x \in P}  =
 \Hom_{RP}(R,-) \colon RP\lmod \to R\lmod$. Writing out the
 definition of $H^0_I$  we see that the composed functor
 $H^0_I \circ \Hom_{RP}(R,-)$ is isomorphic to the functor $\colim_n
 \Hom_{RP}(R/I^n,-)$. Proposition \ref{prop1} also follows from this isomorphism
 and the Grothendieck spectral sequence. On the other hand, since
 filtered colimits commute with finite limits, there is an
 isomorphism $H^0_I \circ \Hom_{RP}(R,-) \cong \Hom_{RP}(R,-) \circ
 H^0_I$. Suppose under the assumptions of Theorem \ref{decomp} that $E$ is a
 projective resolution of $K$ considered as an $KP$-module. The
 Grothendieck spectral sequence for the
 composition $\Hom_{RP}(R,-) \circ
 H^0_I$ then has $E_2$-term
 \begin{eqnarray*}
   E_2^{p,q}(M) &=& (R^p \Hom_{RP}(R,-))((R^qH^0_I)(M)) \\
   &=& H_{-p} \Hom_{RP}(R \otimes_K E,H^q_I(M)) \\
   &\cong& H_{-p}(\bigoplus_{x \in P} \Hom_{KP}(E,K(x,0)) \otimes_K H^q_I(M_x)).
 \end{eqnarray*}
 The statement of Theorem \ref{decomp} implies that this spectral sequence
 collapses at $E_2$.
\end{rem}

\end{document}